\documentclass{jgcc}

\pdfoutput=1
 \usepackage{lastpage}
 \jgccdoi{13}{1}{1}{6521}
 \jgccheading{}{\pageref{LastPage}}{}{}%
 {Jan.~07,~2019}{Feb.~09,~2021}{}

\usepackage{hyperref}
\theoremstyle{plain} 

\newcommand{\Emod}[1]{E(\mathbb{Z}/#1\mathbb{Z})}

\newcommand{\zero}{\mathcal{O}}
\newcommand{\EECN}{Euler elliptic Carmichael number}
\newcommand{\SECN}{strong elliptic Carmichael number}
\newcommand{\exponent}[2]{\epsilon_{#1,#2}(E)}
\newcommand{\exponentNp}{\exponent{N}{p}}
\newcommand{\jacsym}[2]{\left( \frac{#1}{#2} \right)}
\newcommand{\characteristic}{\text{char}}
\newcommand{\tr}{\text{tr}}

\newcommand{\nmid}{\hspace{-4pt}\not|\hspace{2pt}}

\newtheorem*{theorem-non}{Corollary}
\newtheorem*{definition-non}{Definition}

\def\Z{\mathbb{Z}}
\def\Q{\mathbb{Q}}

\def\End{\operatorname{End}}

\def\tr{\operatorname{tr}}
\def\ord{\operatorname{ord}}
\def\Pr{\operatorname{Pr}}
\def\Aut{\operatorname{Aut}}

\let\epsilon\varepsilon

\newif\ifshowold 
\showoldfalse 

\begin{document}

\title[On Types of Elliptic Pseudoprimes]{On Types of Elliptic Pseudoprimes}
\titlecomment{{\lsuper*}{\em AMS Subject Classification}: 14H52, 14K22, 11Y01, 11N25, 11G07, 11G20, 11B99.}

\author[L.~Babinkostova]{Liljana Babinkostova}	
\address{Boise State University}	
\email{liljanababinkostova@boisestate.edu}  

\author[A.~ Hern\'{a}ndez-Espiet]{A. Hern\'{a}ndez-Espiet}	
\address{Rutgers University}	
\email{ah1112@math.rutgers.edu}  

\author[H. Y. ~Kim]{H. Y. Kim}	
\address{University of Wisconsin-Madison}	
\email{hyunjong.kim@math.wisc.edu}  
\thanks{This research was supported by the National Science Foundation under the Grant number DMS-1659872.}	




\begin{abstract}
  \noindent We generalize Silverman's \cite{Korselt} notions of elliptic pseudoprimes and elliptic Carmichael numbers to analogues of Euler-Jacobi and strong pseudoprimes. We inspect the relationships among \EECN{s}, \SECN{s}, products of anomalous primes and elliptic Korselt numbers of Type I, the former two of which we introduce and the latter two of which were introduced by Mazur \cite{Mazur} and Silverman \cite{Korselt} respectively. In particular, we expand upon the work of Babinkostova et al. \cite{REU2016} on the density of certain elliptic Korselt numbers of Type I which are products of anomalous primes, proving a conjecture stated in \cite{REU2016}.
\end{abstract}

\maketitle

\section{Introduction}\label{intro}

The problem of efficiently distinguishing the prime numbers from the composite numbers has been a fundamental problem for a long time. One of the first primality tests in modern number theory came from Fermat Little Theorem: if $p$ is a prime number and $a$ is an integer not divisible by $p$, then $a^{p-1} \equiv 1 \pmod{p}$. The original notion of a pseudoprime (sometimes called a Fermat pseudoprime) involves counterexamples to the converse of this theorem. A {\em pseudoprime} to the base $a$ is a composite number $N$ such $a^{N-1}\equiv 1 \mod N$. A number $N$ which is a pseudoprime to all bases $a$ with $(a, N) = 1$ is a {\em Carmichael number}. Carmichael numbers had been studied by Korselt \cite{KC} who gave the following criterion for Carmichael numbers, which was later rediscovered by Carmichael \cite{C}: A positive composite number $N$ is a Carmichael number if and only if $N$ is odd, square-free, and every prime $p \vert N$ has the property that $(p -1) \vert (N-1)$. In 1986, the long-standing conjecture that there are infinitely many Carmichael numbers was proven by Alford, Granville, and Pomerance \cite {AGP}. In 1989, Gordon introduced the notion of an {\em elliptic pseudoprime} \cite{DG2} as a natural extension of the definition of a pseudoprime from groups arising from elliptic curves with complex multiplication.
\begin{definition-non}[\cite{DG2}]
Let $E/\Q$ be an elliptic curve with complex multiplication by an order in $\Q(\sqrt{-d})$ and let $P \in E(\mathbb{Q})$ has an infinite order. A composite number $N$ is called an \emph{elliptic pseudoprime} if $\jacsym{-d}{N} = -1$, $N$ is coprime to the discriminant of $E$, and $N$ satisfies $(N+1)P\equiv \zero \pmod{N}$. 
\end{definition-non}
As in \cite{Korselt}, we extend the notion of Euler elliptic pseudoprimes and strong Elliptic pseudoprimes, defined by Gordon in \cite{DG2}, to arbitrary elliptic curves $E/\Q$ and points $P \in \Emod{N}$.  
\begin{definition-non}
Let $N \in \Z$, let $E/\Q$ be an elliptic curve, and let $P \in \Emod{N}$. Let$L(E/\Q,s) = \sum_n \frac{a_n}{n^s}$ be the $L$-series of $E/\Q$ and $N+1-a_N$ be even. Then, $N$ is an \emph{Euler elliptic pseudoprime} for $(E,P)$ if $N$ has at least two distinct prime factors, $E$ has good reduction at every prime $p$ dividing $N$, and 
$$ \left( \frac{N+1-a_N}{2} \right) P \equiv \begin{cases} \zero \pmod{N} &\text{if } P = 2Q \text{ for some } Q \in \Emod{N} \\
\text{a 2-torsion point} \pmod N &\text{otherwise}. \end{cases} $$
\end{definition-non}
\begin{definition-non}
Let $E/\Q$ be an elliptic curve with complex multiplication by an order in $\Q(\sqrt{-d})$, let $P$ be a point in $E$ of infinite order, and let $N$ be a composite number with $\gcd(N,6\Delta) = 1$. Let $s$ and $t$ be integers satisfying $N+1 = 2^s  t$, where $t$ is odd. An elliptic pseudoprime $N$ is called a \emph{strong elliptic pseudoprime} for $(E,P)$ if
\begin{enumerate}
\item $t P = \zero \pmod{N}$ or
\item $(2^r t)P$ \text{is a point of order} $2 \pmod N$, for some $r$ with $0\leq r<s$.
\end{enumerate}
\end{definition-non}
We also define the notions of \emph{Euler elliptic Carmichael numbers} and \emph{strong elliptic Carmichael numbers} and identify Korselt criteria for \EECN{s} (Proposition \ref{PropKE}) and \SECN{s} (Proposition \ref{PropKS}). Using these criteria, we show that \SECN{s} are also \EECN{s} when applicable (Corollary \ref{SE}). In Section 4 we investigate the elliptic Korselt numbers of Type I introduced in \cite{Korselt} and show the following result which proves Conjecture 4.9 from \cite{REU2016}.

\begin{theorem-non}
Let $5 \leq p,q \leq M$ be randomly chosen distinct primes and let $N = pq$. Let $\Emod{N}$ be a randomly chosen elliptic curve with good reduction at $p$ and $q$ such that $(p+1-a_p), (q+1-a_q) \mid (N+1-a_N)$. Then
\begin{align*} 
\lim_{M \rightarrow \infty} \Pr[a_p \text{ or } a_q \neq 1] = 0  \text { and } \lim_{M \rightarrow \infty} \Pr[\#\Emod{N} = N+1-a_N] = 1
\end{align*} 

\end{theorem-non}
In Section 5 we investigate the relationship of elliptic Korselt numbers of Type I to \SECN{s} and \EECN{s}. In particular, we show conditions under which elliptic Korselt numbers of Type I are equivalent to strong elliptic Carmichael numbers (Proposition \ref{PropEKIEECN}), as well as conditions under which elliptic Korselt numbers of Type I are equivalent to Euler elliptic Carmichael numbers (Proposition \ref{TypeSECN}).

\section{Preliminaries} 

\subsection{Notation}
For an integer $a$ and a prime $p$, the Legendre symbol $\jacsym{a}{p}$ is defined as
\begin{align*}
	\jacsym{a}{p} = \begin{cases} 0 &\text{if } p \mid a \\1 &\text{if } p \nmid a \text{ and } a \equiv x^2 \pmod{p} \text{ for some } x \in \mathbb{Z}/p\mathbb{Z} \\ -1 &\text{otherwise}. \end{cases}
\end{align*}
For an integer $a$ and a positive odd integer $N$, the Jacobi symbol $\jacsym{a}{N}$ is an extension of the Legendre symbol; if the prime factorization of $N$ is $N = p_1^{e_1} \cdots p_k^{e_k}$, then 
\begin{align*}
	\jacsym{a}{N} = \jacsym{a}{p_1}^{e_1} \cdots \jacsym{a}{p_k}^{e_k}.
\end{align*}
\par
For an integer $N$ and a prime $p$, the $p$-adic order, $\ord_p(N)$, is the largest nonnegative integer $e$ such that $p^e$ divides $N$ if $N \neq 0$ and is $\infty$ otherwise. Given that $e = \ord_p(N)$, we also write $p^e \mid \mid N$. 

\subsection{Elliptic Curves} \label{SectionPrelim}
We introduce some elementary features of elliptic curves which are relevant to the topics presented in this paper. We refer the reader to \cite{Silverman} and \cite{W1} for detailed introduction to elliptic curves. Let $k$ be a field and let $\overline{k}$ denote its algebraic closure. An \emph{elliptic curve} $E$ over a field $k$ is a non-singular {\footnote{\tiny An algebraic curve is said to be non-singular if there is not point on the curve at which all partial derivatives vanish.}} curve with an affine equation of the form
\begin{align} 
	E/k: y^2 + a_1xy + a_3y = x^3 + a_2x^2 + a_4x + a_6     \label{eq1}
\end{align}
where $a_1,a_2,a_3,a_4,a_6 \in k$. An equation of the above form (\ref{eq1}) is called a \emph{generalized Weierstrass} equation. Recall that the points in projective space $\mathbb{P}^2(k)$ correspond to the equivalence classes in $k^3 - \{(0,0,0)\}$ under the equivalence relation $(x, y, z) \sim (ux, uy, uz)$ with $u\in k^{\times}$. The equivalence class containing $(x,y,z)$ is denoted by $[x:y:z]$. The projective equation corresponding to the affine equation (\ref{eq1}) is the homogeneous equation
\begin{align}
	E/k: y^2z + a_1xyz + a_3yz^2 = x^3 + a_2x^2z + a_4xz^2 + a_6z^3,  \label{eq}
\end{align}
where $a_1,a_2,a_3,a_4,a_6 \in \overline{k}$.  
If $\characteristic(k)\neq 2,3 $, then the defining equation of $E$ can be put, after a linear change of variables, in the \emph{Weierstrass normal form}
\begin{align*}
	E/k: y^2 = x^3 + Ax + B  \label{eq2}
\end{align*} 
where $A, B \in k$.
A projective curve $y^2z = x^3+Axz^2+Bz^3$ is non-singular if and only if the discriminant $(4A^3+27B^2) \neq 0$. The points $[x : y : z]$ on the projective curve $$y^2z = x^3+Axz^2+Bz^3$$ are the points [x : y : 1]  where $(x, y)$ is a solution to $y^2 = x^3 + Ax + B$ along with the point $[0 : 1 : 0]$. The point $[0 : 1 : 0]$ is called the \emph{point at infinity} and is denoted by $\mathcal{O}$. The projective points of the elliptic curve $E$ over $\overline{k}$ form an abelian group with $[0 : 1 : 0]$ (point at infinity) as an identity element.  An elliptic curve $E(\mathbb{Z}/N\mathbb{Z})$ is the set of solutions $[x:y:z]$ (insisting that $\gcd(x, y, z, N)=1$) in projective space over $\mathbb{Z}/N\mathbb{Z}$ to a Weierstrass equation $E/k: y^2 + a_1xy + a_3y = x^3 + a_2x^2 + a_4x + a_6$  where the discriminant $\Delta(E)=4A^3+ 27B^2$ has no prime factor in common with $N$. There is a group law on $\Emod{N}$ given by explicit formulae which can be computed (see \cite{W1}). For a given elliptic curve $E/\mathbb{Q} : y^2=x^3 +Ax+B$ where $A, B,N \in \mathbb{Z}$ with $N$ positive odd such that $\gcd(N, \Delta(E))=1$ there is a group homomorphism from $E/\mathbb{Q}$ to $\Emod{N}$ by representing the points in $E/\mathbb{Q}$ as triples $[x:y:z]\in \mathbb{P}^2(k)$.  If the prime factorization of $N$ is $N = p_1^{e_1} \cdots p_k^{e_k}$ then $\Emod{N}$ is isomorphic as a group to the direct product of elliptic curve groups $ \Emod{N} \simeq \Emod{p_1^{e_1}} \oplus \cdots \oplus \Emod{p_k^{e_k}}$. In particular, if we let $E_i$ be the reduction of $E$ modulo $p_i$, then $E_i$ is an elliptic curve over the field $\mathbb{F}_{p_i}$. It can be shown that $\#E(\mathbb{Z}/p_i^{e_i}\mathbb{Z})=p_i^{{e_i}-1}\#E_i(\mathbb{F}_{p_i})$. We refer the reader to \cite{{Lenstra}, {HL}, {W1}} for details about elliptic curves over $\mathbb{Z}/N\mathbb{Z}$.	 Associated to $E/\mathbb{Q}$ is the $L$-function $L(E,s)$, which can be defined as the Euler product
	\begin{align*}
		L(E,s) = \prod_p \frac{1}{1- a_p p^{-s} + 1_E(p) p^{1-2s}}
	\end{align*}
	where
	\begin{align*}
	1_E(p) = \begin{cases} 1 &\text{if } E \text{ has good reduction at } p \\
													0 &\text{otherwise} \end{cases}
\end{align*}
 and $a_p = p+1-\#\Emod{p}$ whether or not $E$ has good reduction at $p$. Alternatively expressing $L(E,s)$ as the Dirichlet series $L(E,s) = \sum_n \frac{a_n}{n^s}$, the map sending a positive integer $n$ to the coefficient $a_n$ is a multiplicative function with
\begin{align*}
	a_{1} &= 1 \\
	a_{p^e} &= a_p a_{p^{e-1}} - 1_E(p) p a_{p^{e-2}} \qquad \text{for all } e \geq 2.
\end{align*}
See \cite[Chapter 8.3]{DiamondShurman} and \cite[Appendix C, Section 16]{Silverman} for more on $L$-series of elliptic curves.  Also, recall that the endomorphism ring $\End(E)$ of $E(\overline{\mathbb{Q}})$ is isomorphic either to $\mathbb{Z}$ or to an order in an imaginary quadratic field, say $\mathbb{Q}(\sqrt{- d})$ where $d$ is a positive squarefree integer. In the latter case, $E$ is said to have complex multiplication in $\mathbb{Q}(\sqrt{-d})$. For a curve with complex multiplication by $\mathbb{Q}(\sqrt{-d})$, $\#E(\mathbb{F}_p)=p+1$ if $p$ does not split in $\mathbb{Q}(\sqrt{-d})$. If $p$ splits in $\mathbb{Q}(\sqrt{-d})$, say $p=\pi\overline{\pi}$, then $\#E(\mathbb{F}_p)=p+1-\tr(u\pi)$ where $u$ is some unit in the field.  In general, $p$ splits in $\mathbb{Q}(\sqrt{-d})$ if $(\frac{-d}{p})=1$. If $E/\mathbb{Q}$ has complex multiplication (CM) by $\mathbb{Q}(\sqrt{-d})$, then $d \in \{1,2,3,7,11,19,43,67,163\}$. Up to isomorphism over $\mathbb{Q}$, there are only thirteen elliptic curves with CM by an order in one of these fields.  

The following fact will be used throughout the paper: let $E/\Q$ be an elliptic curve with complex multiplication in $\Q(\sqrt{-d})$ and let $N > 0$ be an integer such that all its prime factors $p_i >3$ and such that the Jacobi symbol $\left( \frac{-d}{N} \right)= -1$. In this case, there is a prime $p$ such that the $p$-adic order $\ord_p(N)$ is odd and $\left( \frac{-d}{p} \right) = -1$. Thus, $a_p \equiv 0 \pmod{p}$ (see \cite[Proposition 4.31 and Theorem 10.7]{W1}). Moreover, by Hasse's theorem $|a_p| \leq 2\sqrt{p}$ which implies that $a_p = 0$. Since $\ord_p(N)$ is odd, $a_{p^{\ord_p(N)}} = 0$\footnote{\small More generally, $a_{p^{2k+1}} = 0$ and $a_{p^{2k}} = (-p)^k$ for $k \geq 0$ assuming that $a_p = 0$} and since $n \mapsto a_n$ is a multiplicative function, $a_N = 0$. Note that this claim is also true for $p=2, 3$ based on known facts about the l-adic Galois representation attached to $E$ and the trace of the image of Frobenius. We refer the reader to Section 13.2 from \cite{W1} for a detailed  explanation for the above statement for $p=2, 3$.\\

\subsection {Elliptic Pseudoprimes.}  In this section we give some background on elliptic pseudoprimes in general. For other articles that study elliptic pseudoprimes and related notions see \cite{{BM}, {CLS}, {DW}, {Ekstrom}, {EPT}, {GP}, {HL}, {MM}, {Morain}, {Muller}}. Elliptic pseudoprimes are analogous to Fermat pseudoprimes, which are composites $N$ for which $$a^{N-1} \equiv 1 \pmod{N}$$ for a given $a \in \Z/N\Z$. 

In \cite{DG} Gordon introduces the notion of an elliptic pseudoprime as an analog of Fermat pseudoprime.  While the notion of an elliptic pseudoprime in \cite{DG,DG2} is given with respect to an elliptic curve $E/\Q$ and a point $P \in E(\Q)$ of infinite order, we will also apply these definitions to points $P \in \Emod{N}$.
 
\begin{defi}\cite{DG2}\label{GordonEP}
Let $E/\Q$ be an elliptic curve with complex multiplication in $\Q(\sqrt{-d})$, let $P$ be a point in $E$ of infinite order, and let $N$ be a composite number with $\gcd(N,6\Delta) = 1$. Then, $N$ is an \emph{elliptic pseudoprime} for $(E,P)$ if $\jacsym{-d}{N} =-1$ and 
$$(N+1) P \equiv \mathcal{O}\pmod N\footnote{\tiny For details on computing multiples of points in elliptic curve modulo $N$, see \cite[Chapter 3.2]{W1} or Appendix \ref{SectionMult}.}$$
\end{defi}

In \cite{Korselt}, Silverman extends Gordon's aforementioned notion of elliptic pseudoprimes by allowing any elliptic curve $E/\Q$, not just elliptic curves with complex multiplication, as well as any $P \in \Emod{N}$. 

\begin{defi}\cite{Korselt}\label{SilvermanEP}
	Let $N \in \mathbb{Z}$, let $E/\mathbb{Q}$ be an elliptic curve, and let $P \in \Emod{N}$. Write the $L$-series of $E/\mathbb{Q}$ as $L(E/\mathbb{Q},s) = \sum_n \frac{a_n}{n^s}$. Then $N$ is an \emph{elliptic pseudoprime} for $(E,P)$ if $N$ has at least two distinct prime factors, $E$ has good reduction at every prime $p$ dividing $N$, and $(N+1-a_N) P \equiv \zero \pmod{N}$.
\end{defi}
It is not hard to check that for (most) $N$, $\jacsym{-d}{N} = -1$ and $N$ is square-free if and only if $a_N = 0$. Thus, $(n + 1 - a_N)P = (n + 1)P$, so (most) elliptic pseudoprimes in Gordon's sense are also pseudoprimes in Silverman's sense.

Again, $N$ is a pseudoprime in this case because it displays a behavior that it would if it were prime. Indeed, if $N$ is a prime, then $a_N = 0$ as shown in Section \ref{SectionPrelim}. Thus, $\# \Emod{N} = N+1$, so $(p+1)P \equiv \zero \pmod{p}$ for all $P \in \Emod{p}$. $N$ is therefore guaranteed to be composite if $(N+1)P \not\equiv \zero \pmod{N}$, but $N$ may or may not be prime if $(N+1)P \equiv \zero \pmod{N}$.

In \cite{DG2}, Gordon defines also the notion of Euler elliptic pseudoprimes and strong elliptic pseudoprimes, analogously to Euler-Jacobi pseudoprimes and strong pseudoprimes, respectively. Let $p$ be an odd prime and let $a \in \Z/p\Z$ be nonzero. Since $a^{p-1} \equiv 1 \pmod{p}$ and since $\Z/p\Z$ is a field, $a^{\frac{p-1}{2}} \equiv \pm 1 \pmod{p}$. An odd composite integer $N$ is called an {\em Euler pseudoprime} with respect to a nonzero base $a \in \Z/N\Z$ if $a^{\frac{N-1}{2}} \equiv \pm 1 \pmod{N}$. In fact, Euler showed that $a^{\frac{p-1}{2}} \equiv \jacsym{a}{p} \pmod{p}$. This criterion is the basis to the Solovay-Strassen test \cite{SolovayStrassen}. An odd composite integer $N$ is called an {\em Euler-Jacobi pseudoprime} with respect to a nonzero base $a \in \Z/N\Z$ if $a^{\frac{N-1}{2}} \equiv \jacsym{a}{N} \pmod{N}$.  \par
	Strong pseudoprimes are adversaries to the Miller-Rabin primality test \cite{Miller, Rabin}. For an odd prime $p$, express $p-1$ as $p-1 = 2^s t$ where $s,t \in \Z$ with $t$ odd. For any nonzero $a \in \Z/p\Z$, one of the following holds:
	\begin{enumerate}
		\item  $a^t \equiv 1 \pmod{p}$ or
		\item $a^{2^r t} \equiv -1 \pmod{p}$ for some integer $r$ with $0 \leq r < s$.
	\end{enumerate}
	As such, an odd composite number $N$ is a strong pseudoprime for a nonzero base $a \in \Z/p\Z$ if, when expressing $N-1 = 2^s t$ with $t$ odd,
	\begin{enumerate}
		\item $a^t \equiv 1 \pmod{N}$ or 
		\item  $a^{2^r t} \equiv -1 \pmod{N}$ for some integer $r$ with $0 \leq r < s$. 
	\end{enumerate}

Just as in the definition of elliptic pseudoprimes, $N+1$ takes the place of $N-1$ in the definition for Euler elliptic pseudoprime and strong elliptic pseudoprime.

\begin{defi}\cite{DG2}
Let $E/\Q$ be an elliptic curve with complex multiplication in $\Q(\sqrt{-d})$, let $P$ be a point in $E$ of infinite order and let $N$ be a composite number with $\gcd(N,6\Delta) = 1$.
An elliptic pseudoprime $N$ is called an \emph{Euler elliptic pseudoprime} for $(E,P)$ if 
$$ \left( \frac{N+1}{2} \right) P \equiv \begin{cases} \zero \pmod{N} &\text{if } P = 2Q \text{ for some } Q \in \Emod{N} \\																			                          \text{a 2-torsion point modulo } N &\text{otherwise}. \end{cases} $$
\end{defi}

For a prime $p$, recall that the points of order 2 in $\Emod{p}$ are exactly the points of the form $(x,y) = [x:y:1]$ where $2y + a_1x + a_3 \equiv 0 \pmod{p}$. Recall that such points are exactly the points of the form $(x,0) = [x:0:1]$ if $E$ is in Weierstrass normal form. If $P$ is not a double point modulo $N$ and if $\left( \frac{N+1}{2} \right) P$ is not $\zero$ or of the form $[x:y:1]$ where $2y+a_1x+a_3 \equiv 0 \pmod{N}$, then $N$ must be composite. We therefore not consider such an $N$ to be an Euler elliptic pseudoprime, even if $2  \left( \frac{N+1}{2} \right) P \equiv \zero \pmod{N}$. In other words, by a 2-torsion point modulo $N$, we consider the point $\zero$ or a point of the form $[x:y:1]$ where $2y+a_1x+a_3 \equiv 0 \pmod{N}$.
 
For a prime $p$, the points of order 2 in $\Emod{p}$ are exactly the points of the form $(x,0) = [x:0:1]$. If $P$ is not a double modulo $N$ and if $\left( \frac{N+1}{2} \right) P$ is not $\zero$ or of the form $[x:0:1]$, then $N$ must be composite. We will therefore not consider such an $N$ to be an Euler elliptic pseudoprime, even if $2 \left( \left( \frac{N+1}{2} \right) \right) P \equiv \zero \pmod{N}$. In other words, by a 2-torsion point modulo $N$, we will mean $\zero$ or a point of the form $[x:0:1]$. \par

Here, a $2$-torsion point modulo $N$ should not simply be understood as a point $P \in \Emod{N}$ such that $2P \equiv \mathcal{O} \pmod{N}$. Rather, we will say that $P$ is a $2$-torsion point modulo $N$ if (1)$P \equiv \mathcal{O} \pmod{p^e}$ for all $p^e \mid \mid N$, or (2)$P$ has order $2$ modulo $p^e$ for all $p^e \mid \mid N$. If (1) and (2) do not hold, but $2P \equiv \mathcal{O} \pmod{N}$, then a nontrivial factor of $N$ can easily be obtained from the coordinates of $P$.

In \cite{DG2}, Gordon also required that $N\equiv 1 \mod 4$, but in \cite{Muller} it has been shown that this condition is not needed. If $p$ is a prime, then for elliptic curves $E/k: y^2=x^3+Ax+B$ the $2$-torsion points in $E(\mathbb{F}_p)$ (points $P$ such that $2P=\mathcal{O}$) are of the form $(x,0)$, where $x$ is a root of $x^3+Ax+B=0 \mod p$.
In \cite{DG2}, Gordon does not quite define Euler elliptic psuedoprimes as above. If $p$ is a prime and if $\#\Emod{p} = p+1$, then by \cite[Lemma 4.8]{Schoof} we have that $\Emod{p} \simeq \Z/(p+1) \Z$ or $\Z/((p+1)/2) \Z \oplus \Z/2\Z$, with the latter case happening only if $p \equiv 3 \pmod{4}$. Gordon thus puts the additional restriction that $N \equiv 1 \pmod{4}$ and requires that $\left( \frac{N+1}{2} \right) P$ is a 2-torsion point modulo $N$ which is not $\zero$ in the case that $P \not\equiv 2Q$ for all $Q \in \Emod{N}$. Nevertheless, we will allow for $N \equiv 3 \pmod{4}$ when defining Euler elliptic pseudoprimes.

\begin{defi}
Let $E/\Q$ be an elliptic curve with complex multiplication by an order in $\Q(\sqrt{-d})$, let $P$ be a point in $E$ of infinite order, and let $N$ be a composite number with $\gcd(N,6\Delta) = 1$. Further let $s$ and $t$ be integers satisfying $N+1 = 2^s  t$, where $t$ is odd. An elliptic pseudoprime $N$ is called a \emph{strong elliptic pseudoprime} for $(E,P)$ if
\begin{enumerate}
\item $t P = \zero \pmod{N}$ or
\item $(2^r t)P$ is a point of order 2 modulo $N$, for some $r$ with $0\leq r<s$.
\end{enumerate}
\end{defi}
Similarly as before, we will say that a point $P \in \Emod{N}$ is a point of order 2 modulo $N$ if and only if $P$ is of the form $[x:y:1]$ where $2y+a_1x+a_3 \equiv 0 \pmod{N}$. Equivalently, by the Chinese Remainder Theorem, $P$ reduces to a point $[x':y':1]$ modulo $p^e$ such that $2y' + a_1x' + a_3 \equiv 0 \pmod{p^e}$ for every $p^e \mid \mid N$.  \par
Just as before, a point $P \in \Emod{N}$ of order $2$ modulo $N$ is a point such that $P$ has order $2$ modulo $p^e$ for all $p^e \mid \mid N$. 

For Fermat pseudoprimes, all strong pseudoprimes fulfill the corresponding Euler criteria, i.e., are Euler pseudoprimes. However, this doesn't carry over in the case of elliptic pseudoprimes. Just as before, a point $P \in \Emod{N}$ of order $2$ modulo $N$ is a point such that $P$ has order $2$ modulo $p^e$ for all $p^e \mid \mid N$. For Fermat pseudoprimes, all strong pseudoprimes fulfill the corresponding Euler criteria, i.e., are Euler pseudoprimes. However, this doesn't carry over in the case of elliptic pseudoprimes. 

\begin{exa} \label{ExampleMuller}
The following example is a corrected version of the example given in \cite{Muller} and it shows that strong elliptic pseudoprimes do not need to be Euler elliptic pseudoprimes. 
\begin{align*}
 N &= 676258600736819377469073681570025709 \\
	&=	47737 \cdot 275183 \cdot 1212119 \cdot 2489759 \cdot 3178891 \cdot 5366089
\end{align*}
and let $E$ be the curve $E: y^2 = x^3 - 3500 x - 98000$, given in \cite[Table 1]{DG2}, and with complex multiplication in $\Q(\sqrt{-7})$. The example from \cite{Muller} uses the point $P = (84,884)$. However, this point is in fact not in $E$. Rather $P$ should be $(84,448)$.
Note that $N \equiv 1 \pmod{4}$ and $\jacsym{-7}{N} = -1$.  

\begin{align*}
	(N+1)P \equiv \zero \pmod{N},
\end{align*}
so $N$ is an elliptic pseudoprime for $(E,P)$. \par

While the author in \cite{Muller} states that
\begin{align*}
	\left( \frac{N+1}{2} \right) P \equiv (654609963152984637027391710649598749, 0) \pmod{N},
\end{align*}
the point $(654609963152984637027391710649598749, 0)$ is not in the elliptic curve $\Emod{N}$. In fact, 
\begin{align*}
	\left( \frac{N+1}{2} \right) P \equiv (513078336047534294929224848649215641, 0) \pmod{N}.
\end{align*}
Since $\frac{N+1}{2}$ is odd, $N$ is a strong elliptic pseudoprime for $(E,P)$. On the other hand, there is a point 
\begin{align*}
	Q = (427631894156657698513741722706642740,349223536492541846798816891095072158) 
\end{align*}
on $\Emod{N}$ such that 
\begin{align*}
	2Q \equiv (84,448) \equiv P \pmod{N}.
\end{align*}
Thus, $N$ is not an Euler elliptic pseudoprime. For more errors of this kind that appear in \cite{Muller}, see Appendix \ref{SectionMuller}.
\end{exa}

Similarly, the example below shows that Euler elliptic pseudoprimes are not necessarily strong elliptic pseudoprime.

\begin{exa}
Let $N = 7739 = 71 \cdot 109$, $E: y^2 = x^3 -1056 x + 13352$ and $P = (33,121)$. As listed in \cite[Table 1]{DG2}, $E$ has complex multiplication in $\Q(\sqrt{-11})$ and $\jacsym{-11}{N} = -1$. Moreover, $N+1 = 7740 = 2^2 \cdot 1935$. Compute
\begin{align*}
	1935 P &\equiv \zero \pmod{71} \text{ and}\\
	1935 P &\equiv (102,0) \pmod{109},
\end{align*}
so $N$ is not a strong elliptic pseudoprime. However, $N$ is an Euler elliptic pseudoprime because
\begin{align*}
	\left( \frac{N+1}{2} \right) P \equiv (2\cdot 1935) P &\equiv \zero \pmod{N}.
\end{align*}
\end{exa}

\section {Euler elliptic pseudoprimes and Strong elliptic pseudoprimes}
In \cite{Korselt}, Silverman extends Gordon's aforementioned notion of elliptic pseudoprimes by allowing any elliptic curve $E/\Q$, not just elliptic curves with complex multiplication.
\begin{defi}\cite{Korselt}
	Let $N \in \mathbb{Z}$, let $E/\mathbb{Q}$ be an elliptic curve, and let $P \in \Emod{N}$. Write the $L$-series of $E/\mathbb{Q}$ as $L(E/\mathbb{Q},s) = \sum_n \frac{a_n}{n^s}$. Call $N$ an \emph{elliptic pseudoprime} for $(E,P)$ if $N$ has at least two distinct prime factors, $E$ has good reduction at every prime $p$ dividing $N$, and $(N+1-a_N) P \equiv \zero \pmod{N}$.
\end{defi}

As in \cite{Korselt}, we extend Gordon's notions of Euler elliptic pseudoprimes and strong elliptic pseudoprimes, by allowing general elliptic curves over $\Q$ and using $N+1-a_N$ in place of $N+1$. 

\begin{defi}
	Let $N \in \Z$, let $E/\Q$ be an elliptic curve, and let $P \in \Emod{N}$. Write the $L$-series of $E/\Q$ as $L(E/\Q,s) = \sum_n \frac{a_n}{n^s}$ and suppose that $N+1-a_N$ is even. Then, $N$ is an \emph{Euler elliptic pseudoprime} for $(E,P)$ if $N$ has at least two distinct prime factors, $E$ has good reduction at every prime $p$ dividing $N$, and 
	$$ \left( \frac{N+1-a_N}{2} \right) P \equiv \begin{cases} \zero \pmod{N} &\text{if } P = 2Q \text{ for some } Q \in \Emod{N} \\
																															\text{a 2-torsion point modulo } N &\text{otherwise}. \end{cases} $$
\end{defi}

\begin{rem}
Since the definition of Euler elliptic pseudoprime requires the inspection of the multiple $\left( \frac{N+1-a_N}{2} \right) P$, it makes little sense to discuss whether $N$ is an Euler elliptic pseudoprime if $(N+1-a_N)$ is odd.
\end{rem}

\begin{defi}\
	Let $N \in \Z$, let $E/\Q$ be an elliptic curve given by a normal Weierstrass equation, and let $P \in \Emod{N}$. Write the $L$-series of $E/\Q$ as $L(E/\Q,s) = \sum_n \frac{a_n}{n^s}$. Let $s$ and $t$ be integers satisfying $N+1-a_N = 2^s t$, where $t$ is odd. Then, $N$ is a \emph{strong elliptic pseudoprime} for $(E,P)$ if $N$ has at least two distinct prime factors, $E$ has good reduction at every prime $p$ dividing $N$, and 
	\begin{enumerate}
		\item  $tP \equiv \zero \pmod{N}$ or, given that $N+1-a_N$ is even,
		\item $(2^r t)P$ is a point of order 2 modulo $N$ for some $r$ with $0 \leq r < s$.
	\end{enumerate}
\end{defi}
If $(N+1-a_N)$ is odd in the above definition, then condition (ii) above becomes vacuous as $s = 0$. 

Just as Silverman's definition of elliptic pseudoprimes extend Gordon's definition of elliptic pseudoprimes, these definitions of strong and Euler elliptic pseudoprimes extend Gordon's definitions of strong and Euler elliptic pseudoprimes. As such, we can refer to these definitions of elliptic, strong elliptic, and Euler elliptic pseudoprimes without ambiguity. \par

As in \cite{Korselt}, we extend Gordon's notions of Euler elliptic pseudoprimes and strong elliptic pseudoprimes, by allowing general elliptic curves over $\Q$ and using $N+1-a_N$ in place of $N+1$. 
\begin{defi}
Let $N \in \Z$, let $E/\Q$ be an elliptic curve, and let $P \in \Emod{N}$. Write the $L$-series of $E/\Q$ as $L(E/\Q,s) = \sum_n \frac{a_n}{n^s}$ and suppose that $N+1-a_N$ is even. Then, $N$ is an \emph{Euler elliptic pseudoprime} for $(E,P)$ if $N$ has at least two distinct prime factors, $E$ has good reduction at every prime $p$ dividing $N$, and 
$$ \left( \frac{N+1-a_N}{2} \right) P \equiv \begin{cases} \zero \pmod{N} &\text{if } P = 2Q \text{ for some } Q \in \Emod{N} \\
\text{a 2-torsion point modulo } N &\text{otherwise}. \end{cases} $$
\end{defi}

\begin{defi}
Let $N \in \Z$ and let $E/\Q$ be an elliptic curve given by a normal Weierstrass equation, and let $P \in \Emod{N}$. Write the $L$-series of $E/\Q$ as $L(E/\Q,s) = \sum_n \frac{a_n}{n^s}$. Let $s$ and $t$ be integers satisfying $N+1-a_N = 2^s t$, where $t$ is odd. Then, $N$ is a \emph{strong elliptic pseudoprime} for $(E,P)$ if $N$ has at least two distinct prime factors, $E$ has good reduction at every prime $p$ dividing $N$, and 
	\begin{enumerate}
		\item $tP \equiv \zero \pmod{N}$ or, assuming that $N+1-a_N$ is even,
		\item $(2^r t)P$ is a point of order 2 modulo $N$ for some $r$ with $0 \leq r < s$.
	\end{enumerate}
\end{defi}
If $(N+1-a_N)$ is odd in the above definition, then condition (2) above becomes vacuous as $s = 0$.

\begin{defi}\label{StrongCarmichael}
Let $N \in \mathbb{Z}$ and let $E/\Q$ be an elliptic curve. If $N$ is a strong elliptic pseudoprime for $(E,P)$ for every point $P \in \Emod{N}$, then $N$ is a \emph{\SECN} for $E$.
\end{defi}

\section{Korselt Criteria for \EECN{s} and \SECN{s}}
In \cite{Korselt}, Silverman gives a Korselt criterion for elliptic Carmichael numbers. Any number satisfying this elliptic Korselt criterion is an elliptic Carmichael number, but the converse need not be true.
\begin{defi} [\cite{Korselt}]
Let $N \in \mathbb{Z}$, and let $E/\Q$ be an elliptic curve. Then, $N$ is an \emph{elliptic Korselt number for $E$ of Type I} if $N$ has at least two distinct prime factors and, for every prime $p$ dividing $N$,
\begin{enumerate}
	\item $E$ has good reduction at $p$,
	\item $p+1-a_p$ divides $N+1-a_N$, and
	\item $\ord_{p}(a_N-1) \geq \ord_{p}(N) - \begin{cases} 1 &\text{if } a_p \not\equiv 1 \pmod{p} \\
	                                                                                               0 &\text{if } a_p \equiv 1 \pmod{p} \end{cases}$.
\end{enumerate}
\end{defi}
Here, $\ord_p(N)$ denotes the largest nonnegative integer $e$ such that $p^e$ divides $N$ if $N \neq 0$ and $\infty$ otherwise. Assuming that $e = \ord_p(N)$, we also write $p^e \mid \mid N$.  
\begin{prop}[\cite{Korselt}, Proposition 11] \label{PropEKI}
Let $N \in \mathbb{Z}$ be an odd integer and let $E/\Q$ be an elliptic curve. If $N$ is an elliptic Korselt number for $E$ of Type $I$, then $N$ is an elliptic Carmichael number for $E$.
\end{prop}
In \cite{Korselt}, Silverman introduces two notions of elliptic Korselt numbers. Any number satisfying the following elliptic Korselt criterion must be an elliptic Carmichael number, but the converse is not generally true. 
For an integer $N$ and a prime $p$, the $p$-adic order, $\ord_p(N)$, is the largest nonnegative integer $e$ such that $p^e$ divides $N$ if $N \neq 0$ and is $\infty$ otherwise. Given that $e = \ord_p(N)$, we also write $p^e \mid \mid N$. 

Silverman's second elliptic Korselt criterion gives a necessary and sufficient condition for an integer to be an elliptic Carmichael number for an elliptic curve. In doing so, we will use the following notation, as he does in \cite[Page 8]{Korselt}, for the exponent of a group:
\begin{defi}
	For a group $G$, denote $\epsilon(G)$ as the \emph{exponent of $G$}, i.e. the least positive integer such that $g^{\epsilon(G)} = 1$ for all $g \in G$. Equivalently, $\epsilon(G)$ is the least common multiple of the orders of all of the elements of $G$. \par
	For an elliptic curve $E/\Q$, an integer $N$, and a prime $p$ dividing $N$ at which $E$ has good reduction, write
	$$\exponentNp = \epsilon \left( E\left( \Z/p^{\ord_{p}(N)} \Z \right) \right). $$
\end{defi}

\begin{defi}
	Let $N \in \mathbb{Z}$ and let $E/\Q$ be an elliptic curve. We say that $N$ is an \emph{elliptic Korselt number for E of type II} if $N$ has at least two distinct prime factors and if, for every prime $p$ dividing $N$, 
	\begin{enumerate}
		\item  $E$ has good reduction at $p$ and
		\item $\exponentNp$ divides $N+1-a_N$.
	\end{enumerate}
\end{defi}

\begin{prop}[\cite{Korselt}, Proposition 12]
	Let $N > 2$ be an odd integer, and let $E/\Q$ be an elliptic curve. Then, $N$ is an elliptic Carmichael number for $E$ if and only if $N$ is an elliptic Korselt number for $E$ of type $II$.
\end{prop}

We introduce the notion of \EECN{s} and \SECN{s} and for each of them we give necessary and sufficient Korselt criteria (Proposition \ref{PropKE} and Proposition \ref{PropKS} respectively), akin to the Korselt criterion.
\begin{defi}\label{EulerCarmichael}
Let $N \in \mathbb{Z}$ and let $E/\Q$ be an elliptic curve. If $N$ is an Euler elliptic pseudoprime for $(E,P)$ for every point $P \in \Emod{N}$, then $N$ is an \emph{\EECN} for $E$.
\end{defi}
\begin{defi}\label{StrongCarmichael}
Let $N \in \mathbb{Z}$ and let $E/\Q$ be an elliptic curve. If $N$ is a strong elliptic pseudoprime for $(E,P)$ for every point $P \in \Emod{N}$, then $N$ is a \emph{\SECN} for $E$.
\end{defi}

\begin{prop} \label{PropKE}
Let $N \in \Z$ be an integer with at least two distinct prime factors, let $E/\Q$ be an elliptic curve, and suppose that $N+1-a_N$ is even. Then, $N$ is {an} \EECN{} if and only if $E$ has good reduction at $p$ for every for every prime $p\mid N$ and $\exponentNp\mid \frac{N+1-a_N}{2}$.	
\end{prop}
\begin{proof}
Suppose that $E$ has good reduction at $p$ and $\exponentNp\mid\frac{N+1-a_N}{2}$ for all prime powers $p^e \mid \mid N$. For all $P \in \Emod{N}$, $\left( \frac{N+1-a_N}{2} \right) P \equiv \zero \pmod{p^e}$, so $\left( \frac{N+1-a_N}{2} \right) P \equiv \zero \pmod{N}$. Conversely, suppose that $N$ is {an} \EECN{} for $E$. In particular, $E$ has good reduction at every prime dividing $N$. For each prime power $p^e \mid \mid N$, there is an element of $\Emod{p^e}$ of order $\exponentNp$. Let $P$ be a point of $\Emod{N}$ such that $P$ has order $\exponentNp$ modulo $p^e$ for all prime powers $p^e \mid \mid N$. If $\exponentNp$ is odd for every prime $p$ dividing $N$, then $P \equiv 2Q \pmod{N}$ for some $Q \in \Emod{N}$. Therefore, $\left( \frac{N+1-a_N}{2} \right) P \equiv \zero \pmod{N}$, so $\exponentNp$ must divide $\frac{N+1-a_N}{2}$ for all primes $p$ dividing $N$. Next assume that there are prime powers $p^e \mid \mid N$ such that $\exponentNp$ is even. In this case, $P$ is not a double modulo $p^e$ whenever $\exponentNp$ is even, so $P$ is not a double modulo $N$. Since $N$ is {an} \EECN{} for $E$, $\left( \frac{N+1-a_N}{2} \right) P$ is a 2-torsion point modulo $N$. If $\left( \frac{N+1-a_N}{2} \right) P \equiv \zero \pmod{N}$, then $\exponentNp \mid \frac{N+1-a_N}{2}$ for all primes $p$ dividing $N$, which is the desired result.  Suppose for contradiction that $\left( \frac{N+1-a_N}{2} \right) P$ has order 2 modulo $N$. Let $P'$ be a point of $\Emod{N}$ which satisfies
	\begin{align*}
		P' \equiv \begin{cases} 2P \pmod{p^e} &\text{if } p^e \mid \mid N \text{ with } \exponentNp \text{ even} \\
														\phantom{2}P \pmod{p^e} &\text{if } p^e \mid \mid N \text{ with } \exponentNp \text{ odd}. \end{cases}
	\end{align*}
Note that $P'$ is a double modulo $p^e$ for every prime power $p^e \mid \mid N$ as all points of $\Emod{p^e}$ are doubles if $\exponentNp$ is odd. Therefore, $\left( \frac{N+1-a_N}{2} \right) P' \equiv \zero \pmod{N}$, but 
	\begin{align*}
	\left( \frac{N+1-a_N}{2} \right) P' \equiv \left( \frac{N+1-a_N}{2} \right) P \not\equiv \zero \pmod{p^e}
	\end{align*}
for every prime power $p^e \mid \mid N$ such that $\exponentNp$ is odd. There is thus no prime $p$ dividing $N$ for which $\exponentNp$ is odd. 
Fix a prime power $p_1^{e_1} \mid \mid N$. Now let $P'$ be a point of $\Emod{N}$ which satisfies
	\begin{align*}
		P' \equiv \begin{cases} 2P \pmod{p^e} &\text{if } p = p_1, e = e_1 \\
														\phantom{2}P \pmod{p^e} &\text{if } p^e \mid \mid N \text{ with } p \neq p_1. \end{cases}
	\end{align*}
	Since $N$ has at least two distinct prime factors and $\exponentNp$ is even for all primes $p$ dividing $N$, $P'$ is not a double in $\Emod{N}$. Therefore, $\left( \frac{N+1-a_N}{2} \right) P'$ is a 2-torsion point. However,  
	\begin{align*}
	\left( \frac{N+1-a_N}{2} \right) P' \equiv 2 \left( \left( \frac{N+1-a_N}{2} \right) P\right) \equiv \zero \pmod{p_1^{e_1}}, 
	\end{align*} 
However,
	\begin{align*}
	\left( \frac{N+1-a_N}{2} \right) P' \equiv \left( \left( \frac{N+1-a_N}{2} \right) P\right) \not\equiv \zero \pmod{p^e}
	\end{align*} 
for all prime powers $p^e \mid \mid N$ different from $p_1^{e_1}$, which is a contradiction. Hence, $\left( \frac{N+1-a_N}{2} \right) P$ does not have order 2 modulo $N$, i.e. $\exponentNp \mid \frac{N+1-a_N}{2}$ for all primes $p$ dividing $N$.
\end{proof}

\begin{prop} \label{PropKS}
Let $N \in \Z$ be an odd integer with at least two distinct prime factors, let $E/\Q$ be an elliptic curve, and let $s$ and $t$ be integers satisfying $N+1-a_N = 2^s t$ where $t$ is odd. Then, $N$ is a \SECN{} if and only if, for every prime $p$ dividing $N$, $E$ has good reduction at $p$ and $\exponentNp$ divides $t$. 
\end{prop}
\begin{proof}
Suppose that $E$ has good reduction at $p$ and that $\exponentNp$ divides $t$ for all $p^e \mid\mid N$. Since $\exponentNp$ is the exponent of $\Emod{p^{\ord_p(N)}}$, $tP \equiv \zero \pmod{p^e}$ for every $P \in \Emod{N}$. By the Chinese Remainder Theorem, $tP \equiv \zero \pmod{N}$, so $N$ is a \SECN{}. Conversely, suppose that $N$ is a \SECN{} for $E$. In particular, $E$ has good reduction at every prime dividing $N$. There is an element of $\Emod{p^{\ord_p(N)}}$ of order $\exponentNp$. Let $P$ be a point of $\Emod{N}$ such that $P$ has order $\exponentNp$ modulo $p^e$ for all $p^e \mid\mid N$.  Suppose, towards a contradiction, that $\exponentNp \nmid t$ for some prime $p\mid N$. Consequently, $tP \not\equiv \zero \pmod{N}$. Since $N$ must be a strong elliptic pseudoprime for $(E,P)$, there is $r\in \mathbb{Z}$ satisfying $0 \leq r < s$ for which $(2^r t) P$ is a point of order 2 modulo $N$. Also, there is $p^{e} \mid \mid N$ such that $tP \not\equiv \zero \pmod{p^{e}}$. In fact, this must hold for all $p^e \mid \mid N$; otherwise, $(2^r t) P \equiv \zero \pmod{p^e}$, so $(2^r t) P$ would not be a point of order 2 modulo $p^e$.
Choose  $p_1^{e_1} \mid \mid N$. Let $P'$ be a point of $\Emod{N}$ which satisfies
	\begin{align*}
		P' \equiv \begin{cases} 2P \pmod{p^e} &\text{if } p = p_1, e = e_1 \\
		P \phantom{2} \pmod{p^e} &\text{if } p^e \mid\mid N \text{ with } p \neq p_1. \end{cases}
	\end{align*}
Note that $tP'\equiv \zero \pmod p^e$ for all $p^e \mid \mid {N}$ with $p \neq p_1$. We show that there is no integer $r'$ satisfying $0 \leq r' < s$ for which $(2^{r'}t)P'$ is a point of order 2 modulo $p^e$ for all $p^e \mid \mid N$. In the case when $r' = r$, we have $(2^rt)P' \equiv \zero \pmod {p_1^{e_1}}$ and is of order 2 modulo $p^e$ for all $p^e \mid \mid N$ with $p \neq p_1$. If $r' > r$, then $(2^{r'}t) P' \equiv \zero \pmod{N}$. On the other hand, if $r' < r$, then $(2^{r'}t) P'$ has order greater than 2 modulo $p^e$ for all $p^e \mid \mid N$ with $p \neq p_1$. Thus, there is no such $r'$ and so $N$ is not a strong elliptic pseudoprime for $(E,P')$, which is a contradiction. Hence, $\exponentNp\vert t$ for all primes $p\vert N$ that divide $N$.
\end{proof}

 \begin{rem}
Let $N$ be a composite number which is either not  \EECN{} or not a \SECN{}. In the above propositions, we guarantee the existence of a point $P \in \Emod{N}$ for which $N$ is not an Euler elliptic number (strong elliptic number) for $(E,P)$. However, this does not guarantee the existence of a point $P \in \Emod{N}$ for which $N$ is not an Euler elliptic number (strong elliptic number) for $(E,P)$ and such that $P \not\equiv \zero \pmod{p^e}$ for all prime powers $p^e \mid \mid N$.  On the other hand, if $\exponentNp > 2$ for all primes $p\mid N$, then there is a point $P \in E(\mathbb{Z}/N\mathbb{Z})$ such that $P \not \equiv \zero \pmod{p^e}$ for all prime powers $p^e \mid \mid N$. With $P$ and $P'$ defined to be points of $\Emod{N}$ as in the proofs of Proposition \ref{PropKE} and Proposition \ref{PropKS}, we have $ P,P' \not\equiv \zero \pmod{p^e}$ for all prime powers $p^e \mid \mid N$. For a prime $p \geq 11$, we show that $\epsilon\left( \Emod{p} \right) > 2$. By Hasse's Theorem, $\#\Emod{p} \geq p+1-2\sqrt{p} = (\sqrt{p}-1)^2 >  4$. Therefore, $\#\Emod{p}$ must be divisible either by an odd prime or be a power of $2$ in which case $\#\Emod{p}>4$. Since $\Emod{p}$ is generated by at most 2 elements, the exponent of $\Emod{p}$, that is, $\epsilon\left( \Emod{p} \right)>2$ such that $N$ is not an Euler elliptic number (strong elliptic number) for $(E, P)$
 \end{rem}

To show the existence of \SECN{s} we first define the notion of anomalous primes, introduced by Mazur \cite{Mazur}.
\begin{defi} [\cite{Mazur}]
Let $E/\Q$ be an elliptic curve and let $p$ be a prime number at which $E$ has good reduction. If $p$ divides $\#\Emod{p}$, then $p$ is said to be an \emph{anomalous prime}. By Hasse's theorem, when $p \geq 5$, this is equivalent to saying that $\#\Emod{p} = p$.	
\end{defi}
The proofs of the next two statements are straightforward and are omitted.
\begin{cor}
Let $E/\Q$ be an elliptic curve and let $N = p_1 \cdots p_k$ where $p_1,\ldots,p_k > 3$ are distinct anomalous primes for $E$. Then, $N$ is a \SECN{} for $E$.
\end{cor}
\begin{cor}\label{SE}
Let $E/\Q$ be an elliptic curve and let $N$ be a \SECN{}. If $N+1-a_N$ is even, then $N$ is also an \EECN{}.
\end{cor}
Note that \EECN{s} and \SECN{s} behave differently under Gordon's definition of elliptic pseudoprimes (Def. \ref{GordonEP}).
\begin{exa}
There exist \EECN{s} under Gordon's conditions (Definition \ref{GordonEP}), i.e. that $E$ has complex multiplication in $\Q(\sqrt{-d})$, $\gcd(N,6\Delta) = 1$, and $\left( \frac{-d}{N} \right) = -1$. Let $E: y^2 = x^3 + 80$ be an elliptic curve with complex multiplication in $\Q(\sqrt{-3})$ and let $N = 6119 = 29 \cdot 211$. We have that $\jacsym{-d}{N} = -1$, $\exponent{N}{29} = 30$ and $\exponent{N}{211} = 15$. Moreover, since $\frac{N+1-a_N}{2} = 3060$, $\exponent{N}{p} \mid \frac{N+1-a_N}{2}$ for $p = 29$ and for $p= 211$. 
\end{exa}
\begin{cor}
Let $E/\Q$ be an elliptic curve with complex multiplication in $\Q(\sqrt{-d})$, let $N$ be a composite number with $\gcd(N,6\Delta) = 1$ and $\left( \frac{-d}{N} \right) = -1$. Then, $N$ is not a \SECN.
\end{cor}
\begin{proof}
Since $\left( \frac{-d}{N} \right) = -1$, there is some prime $p$ dividing $N$ for which $\left( \frac{-d}{p} \right) = -1$. In particular, $a_p = 0$, so $\#\Emod{p} = p+1$. The exponent $\exponentNp$ of $\Emod{p^{\text{ord}_p(N)}}$ is therefore even, which implies that $\exponentNp \nmid t$ as $t$ is odd.
\end{proof}

\section{Relationship between \EECN{s}, \SECN{s} and Elliptic Korselt numbers of Type I}
By Proposition \ref{PropEKI}, elliptic Korselt numbers for $E/\Q$ of Type I are elliptic Carmichael numbers, but elliptic Carmichael numbers are generally not elliptic Korselt numbers for $E/\Q$ of Type I. The same holds true for \EECN{s} and \SECN{s}, so we consider the relationships of \EECN{s} and \SECN{s} to elliptic Korselt numbers of Type I. 
\begin{exa}
As in \cite[Example 19]{Korselt}, let $E$ be the elliptic curve $E: y^2 = x^3 + 7x + 3$ and $N = 27563 = 43 \cdot 641$, which is a Type I Korselt number for $E$. We have $a_{43} = 2$, $a_{641} = -15$, $\exponent{N}{43} = 42$ and $\exponent{N}{657} = 657$, so $a_N = -30$. Note that $\left( \frac{N+1-a_N}{2} \right) = 13797$, but $42$ does not divide $13797$. Therefore, $N$ is neither {an} \EECN{} nor a \SECN{} for $E$. 
\end{exa}
\begin{prop} \label{PropEKIEECN}
Let $E/\Q$ be an elliptic curve and let $N$ be an elliptic Korselt number of Type $I$ for $E$. Suppose that $N+1-a_N$ is even. Then, $N$ is {an} \EECN{} for $E$ if and only if, for every prime $p$ dividing $N$, 	
	\begin{enumerate}
		\item  $\left( p+1-a_p \right) \mid \left( \frac{N+1-a_N}{2} \right)$ or
		\item $\Emod{p}$ has exactly three elements of order $2$. 
	\end{enumerate} 	
\begin{proof}
For a fixed prime $p$ dividing $N$, express the cyclic group decomposition of $\Emod{p}$ as $\Emod{p} \simeq \Z/\delta\Z \oplus \Z/\epsilon\Z$ where $\delta \mid \epsilon$. In particular, $p+1-a_p = \#\Emod{p} = \delta\epsilon$ and $\epsilon$ is the exponent of $\Emod{p}$. \par
Suppose that $N$ is an elliptic Korselt number of Type I and {an} \EECN{} for $E$. Let $p$ be a prime dividing $N$ and further suppose that $\left( p+1-a_p \right) \nmid \left( \frac{N+1-a_N}{2} \right)$. We show that $\Emod{p}$ has exactly three elements of order $2$. Since $N$ is an elliptic Korselt number of Type I for $E$, $\left( p+1-a_p \right) \mid (N+1-a_N)$. Therefore, $\ord_{2}(p+1-a_p) = \ord_{2}(N+1-a_N)$. \par
Suppose for contradiction that $p+1-a_p \equiv 0 \pmod{p}$, i.e. $a_p \equiv 1 \pmod{p}$. If $a_p = 1$, then $\#\Emod{p} = p+1-a_p = p$, so $\epsilon = p$. Since $N$ is odd and since $p\vert (N+1-a_N)$, $p$ must divide $\left( \frac{N+1-a_N}{2} \right)$, which is a contradiction. Thus, $a_p \neq 1$. 

If $p \geq 7$, then $a_p \equiv 1 \pmod{p}$ is equivalent to $a_p = 1$ as $|a_p| \leq 2\sqrt{p}$ by Hasse's Theorem, and so $p \leq 5$.  It is easy to see that $\#\Emod{p} = p+1-a_p = 2p$. On the other hand, $\#\Emod{p} = \delta \epsilon$ and $\delta \mid \epsilon$, and so $\delta = 1$ and $\epsilon = 2p$. In particular, $\epsilon = p+1-a_p$. Recall that $\exponentNp$ is the exponent of $\Emod{p^{\ord_p(N)}}$, and so $\epsilon$ divides $\exponentNp$. Since $N$ is {an} \EECN{} for $E$, $\exponentNp \mid \left( \frac{N+1-a_N}{2} \right)$. However, $\epsilon = p+1-a_p$, and so $(p+1-a_p) \mid \left( \frac{N+1-a_N}{2} \right)$, which is a contradiction. Hence, $p+1-a_p \not\equiv 0\pmod{p}$, and so $p+1-a_p$ is not divisible by $p$. \par
Next, suppose that $\delta$ is odd. Since $\delta \epsilon = p+1-a_p$, $\ord_2(\epsilon) = \ord_2(p+1-a_p)$. Moreover, by \cite[The discussion leading up to Proposition 16]{Silverman}, $\exponentNp = p^e \epsilon$ for some nonnegative integer $e$.  In particular, $\ord_2(\epsilon) = \ord_2(\exponentNp)$. Since $\exponentNp \mid \left( \frac{N+1-a_N}{2} \right)$, $\ord_2(p+1-a_p) = \ord_2(\epsilon) = \ord(\exponentNp) < \ord_2(N+1-a_N)$, which contradicts that $\ord_2(p+1-a_p) = \ord_2(N+1-a_N)$. Hence, $\delta$ is even. 
Since $\delta$ is even and $\delta$ divides $\epsilon$, $\epsilon$ must be even. In particular, the $2$-torsion subgroup of $\Emod{p}$ is isomorphic to $\Z/2\Z \oplus \Z/2\Z$. Therefore, there are exactly three points of order $2$ in $\Emod{p}$ as desired. 

Conversely, suppose that $N$ is an elliptic Korselt number of Type I such that \emph{(i)} or \emph{(ii)} holds for every prime $p$ dividing $N$. Since $N$ is an elliptic Korselt number of Type I, an argument in \cite[Equations (4.4) and (4.6)]{Korselt} shows that $p^{\ord_p(N)-1} (p+1-a_p) \mid (N+1-a_N)$. \cite[Remark 14]{Korselt} further gives an exact sequence
	\begin{equation} \label{EqES}
	\begin{aligned}
		0 \rightarrow p\Z/p^{\ord_p(N)} \Z \rightarrow \Emod{p^{\ord_p(N)}} \rightarrow \Emod{p} \rightarrow 0.
	\end{aligned}
	\end{equation}
Suppose that $p+1-a_p$ is not divisible by $p$. In this case, $\Emod{p^{\ord_p(N)}} \simeq \Z/p^{\ord_{p}(N)-1}\Z \oplus \Emod{p}$, and so $\exponentNp = p^{\ord_{p}(N)-1} \epsilon$, where $\epsilon$ is the exponent of $\Emod{p}$ as before. \par	
We will show that $\epsilon \mid \left( \frac{N+1-a_N}{2} \right)$. If $(p+1-a_p) \mid \left( \frac{N+1-a_N}{2} \right)$, then $\epsilon \mid \left( \frac{N+1-a_N}{2} \right)$ because $\epsilon \mid (p+1-a_p)$. On the other hand, if $\Emod{p}$ has exactly three elements of order $2$, then the $2$-torsion subgroup of $\Emod{p}$ is isomorphic to $\Z/2\Z \oplus \Z/2\Z$. In particular, $\delta$ is even. Since $\delta\epsilon = \#\Emod{p} = p+1-a_p$, $\epsilon$ divides $\left( \frac{p+1-a_p}{2} \right)$. Either way, $\epsilon$ divides $\left( \frac{N+1-a_N}{2} \right)$ as desired. Recall that $p^{\ord_{p}(N)-1} \mid (N+1-a_N)$, so $p^{\ord_{p}(N)-1} \mid \left( \frac{N+1-a_N}{2} \right)$. Therefore, $\exponentNp \mid \left( \frac{N+1-a_N}{2} \right)$. Now suppose that $p+1-a_p$ is divisible by $p$. By \cite[Proposition 16]{Korselt}, $p+1-a_p = p$ or $2p$. Since $\delta \mid \epsilon$ and $\delta\epsilon = p+1-a_p$, $\delta = 1$ and $\epsilon = p+1-a_p$. Therefore, $\Emod{p}$ does not have exactly three elements of order $2$, so $(p+1-a_p) \mid \left( \frac{N+1-a_N}{2} \right)$. Recall that $p^{\ord_p(N)-1} (p+1-a_p) \mid (N+1-a_N)$ and since $p$ is odd, $p^{\ord_p(N)-1} (p+1-a_p) \mid \left( \frac{N+1-a_N}{2} \right)$. \cite[Proposition 16]{Korselt} shows that $\exponentNp \mid p^{\ord_p(N)-1} (p+1-a_p)$, so $\exponentNp \mid \left( \frac{N+1-a_N}{2} \right)$ as desired.
\end{proof}
\end{prop}
The following statement summarizes when elliptic Korselt numbers of Type I are strong elliptic Carmichael numbers.
\begin{cor} \label{TypeSECN}
Let $E/\Q$ be an elliptic curve and let $N$ be an elliptic Korselt number of Type $I$ for $E$. Then, $N$ is a \SECN{} for $E$ if and only if $p+1-a_p$ is odd for all primes $p$ dividing $N$.
\begin{proof}
If $p+1-a_p$ is odd for all primes $p$ dividing $N$, then $\exponentNp$ is also odd because $\exponentNp \mid (p+1-a_p)$. Moreover, $(p+1-a_p) \mid (N+1-a_N)$ because $N$ is an elliptic Korselt number of Type I for $E$, so $\exponentNp$ divides the largest odd factor of $N+1-a_N$. By Proposition \ref{PropKS}, $N$ is a \SECN{} for $E$. 
If $p+1-a_p$ is even for some prime $p \mid N$, then $\exponentNp$ must be even. Therefore, $\exponentNp$ does not divide the largest odd factor of $N+1-a_N$, and so $N$ is not a \SECN{} for $E$ by Proposition \ref{PropKS}.
\end{proof}
\end{cor}

\section{Properties of Elliptic Korselt Numbers of Type I}

In \cite[Proposition 4.3]{REU2016} the authors show that products of distinct anomalous primes for an elliptic curve $E/\Q$ are elliptic Korselt numbers of Type I for $E$. Here, we consider the question of how often an elliptic Korselt number of Type I is a product of distinct anomalous primes. We prove the following conjecture from \cite{REU2016}, which deals with the case in which the number in question is semiprime. 

\begin{conj}\label{Conjecture}
For $M \geq 7$, let $5 \leq p,q \leq M$ be distinct primes chosen uniformly at random, and let $N = pq$. Let $\Emod{N}$ be an elliptic curve with good reduction at $p$ and $q$ chosen uniformly at random and $\#\Emod{p} = p+1-a_p$ and $\#\Emod{q} = q+1-a_q$ both divide $N+1-a_N$. Then
\begin{align*}
	\lim_{M \rightarrow \infty} \text{Pr} [\# \Emod{N} = N+1-a_N] = 1.
\end{align*}
\end{conj}
Note that $N = pq$ is an elliptic Korselt number of Type I if and only if $\#\Emod{p}$ and $\#\Emod{q}$ divide $N+1-a_N$ by \cite[Proposition 4.11]{REU2016}. 

\subsection{Bounds on the number of elliptic curves modulo $p$ of prescribed order}
We use Deuring's theorem \cite{D} (see also \cite{Lenstra}), to obtain bounds on the number of elliptic curves modulo $p$ having prescribed order. Write a nonzero integer $\Delta$ as $\Delta = \Delta_0 f^2$ where $\Delta_0$ is square free.
Let $L \left(s, \jacsym{\cdot}{|\Delta_0|} \right)$ be the L-function
\begin{align*}
	L\left(s, \jacsym{\cdot}{|\Delta_0|} \right) = \sum_{n=1}^\infty \frac{\jacsym{n}{|\Delta_0|}}{n^s}
\end{align*}
and let $\psi(f)$ be the multiplicative function defined by
\begin{align*}
	\psi(p^k) = \begin{cases} \frac{p-p^{-k}}{p-1} &\text{if } \jacsym{p}{|\Delta_0|} = 0 \\
														1 &\text{if } \jacsym{p}{|\Delta_0|} = 1\\
														\frac{p+1-2p^{-k}}{p-1} &\text{if } \jacsym{p}{|\Delta_0|} = -1 \end{cases}.
\end{align*}
Recall that the Kronecker class number is $H(\Delta) = \frac{\sqrt{|\Delta|}}{2\pi} L\left(1, \jacsym{\cdot}{|\Delta_0|} \right) \psi(f)$.

\begin{lem} \label{LemmaNumAp}
The number of isomorphism classes of elliptic curves $\Emod{p}$ such that $\#\Emod{p} = p+1-t$ is $H(t^2-4p)$.
\end{lem}
We will use upper and lower bounds for $H(\Delta)$ to prove Conjecture \ref{Conjecture}. Using \cite[Theorem 328]{HardyWright}, one can show that
\begin{align*}
	1 \leq \psi(f) \leq \left( \frac{f}{\varphi(f)} \right)^2 = O\left( (\log \log f)^2 \right).
\end{align*}
where $\varphi$ is the  totient function. 
Since $\Delta_0$ is square free, $\jacsym{\cdot}{|\Delta_0|}$ is a primitive Dirichlet character. The following is a classical result on the upper bound of $L\left(1, \jacsym{\cdot}{|\Delta_0|} \right)$.
\begin{lem}
	$L\left(1, \jacsym{\cdot}{|\Delta_0|} \right) = O(\log |\Delta_0| )$.
\end{lem}	
Moreover, Siegel's Theorem \cite{Siegel} \footnote{\tiny See \cite[Chapter 21]{Davenport}} yields that
\begin{align*}
L\left(1, \jacsym{\cdot}{|\Delta_0|} \right) = \Omega \left( \frac{1}{|\Delta_0|^{\epsilon}} \right)
\end{align*}
for every $\epsilon > 0$. Assuming the generalized Riemann hypothesis, this result can be strengthened as
\begin{align*}
L\left(1, \jacsym{\cdot}{|\Delta_0|} \right) = \Omega\left( \frac{1}{\log\log |\Delta_0|} \right).
\end{align*}

\begin{lem} \label{LemmaBounds}
For all $\epsilon > 0$, 
	\begin{align*}
	|\Delta|^{1/2-\epsilon} \ll H(\Delta) \ll \Delta^{1/2} \log |\Delta| (\log\log |\Delta|)^2
	\end{align*}
In particular, for all $\epsilon > 0$, 
	\begin{align*}
	 |\Delta|^{1/2-\epsilon} \ll H(\Delta) \ll |\Delta|^{1/2+\epsilon}
	 \end{align*}
\end{lem}

\begin{cor} \label{corProb}
Let $N = pq$ be a product of distinct primes and let  $|a_p| \leq 2\sqrt{p}$ and $|a_q| \leq 2\sqrt{q}$. The probability that a randomly chosen elliptic curve $\Emod{N}$ satisfies $\#\Emod{p} = p+1-a_p$ and $\#\Emod{q} = q+1-a_q$ is 
\begin{align*}
	O\left( \frac{ (4p - a_p^2)^{1/2+\epsilon} (4q - a_q^2)^{1/2+\epsilon}}{pq} \right)
\end{align*}
and
\begin{align*}
	\Omega\left( \frac{ (4p - a_p^2)^{1/2-\epsilon} (4q - a_q^2)^{1/2-\epsilon}}{pq} \right)
\end{align*}
for all $\epsilon > 0$. In particular, the probability is
\begin{align*}
	O\left( \frac{(4q - a_q^2)^{1/2-\epsilon}}{p^{1/2-\epsilon} q} \right)
\end{align*}
and
\begin{align*}
	O\left( \frac{ 1}{(pq)^{1/2-\epsilon}} \right).
\end{align*}
\begin{proof}
For a prime $p$, the number of automorphisms on an elliptic curve $\Emod{p}$ is bounded above by $6$. Furthermore, the number of elliptic curves in an isomorphism class with representative $E$ is $(p-1)/\#\Aut E$. Thus there are  $\Theta(p)$ elliptic curves in each isomorphism class. Also, there $p^2-p$ elliptic curves  $\Emod{p}$  with good reduction at $p$. By the Chinese Remainder Theorem, there are $\Theta(p^2q^2)$ elliptic curves $\Emod{N}$  with good reduction at $p$ and $q$. By Lemma \ref{LemmaNumAp}, the number of isomorphism classes of elliptic curves with order $p+1-a_p$ is $H(4p - a_p^2)$. The desired result holds by Lemma \ref{LemmaBounds}.
\end{proof}
\end{cor}

\subsection{The proportion of choices for $p,q,E$ such that $p$ and $q$ are anomalous primes for $E$}
We compute the probability that $p$ and $q$ are anomalous for $E$ assuming that $p$ and $q$ are random distinct primes $5 \leq p,q \leq M$ and assuming that $\Emod{N}$ is any random curve. We will show that
\begin{align*}
	\Pr[a_p \text{ or } a_q \neq 1 \text{ and } (p+1-a_p), (q+1-a_q) \mid (N+1-a_N) ] = o(\Pr[a_p, a_q = 1])
\end{align*}
with respect to $M$. 
\begin{lem} \label{LemmaLBound}
Let $N = pq$ be a product of distinct primes $5 \leq p,q \leq M$ chosen at random. Let $\Emod{N}$ be an elliptic curve with good reduction at $p$ and $q$. The probability that $a_p = a_q = 1$ is 
	$ \Omega\left( \frac{1}{M^{1+\epsilon}} \right)$
for all $\epsilon > 0$.  
\end{lem}
\begin{proof}
The number of primes below $M$ is approximately $\frac{M}{\log M}$. Thus, the number of possible pairs of distinct $p$ and $q$ is $\Theta \left( \frac{M^2}{\log ^{2}M} \right)$, and so
	\begin{align*}
		\Pr[p = p_0, q = q_0 \text{ and } a_p = a_q = 1] = \Omega \left( \frac{1}{p^{1/2 + \epsilon} q^{1/2 + \epsilon} M^2} \right).
	\end{align*}
Then
	\begin{equation} \label{EqProb1}
	\begin{aligned}
		\Pr[a_p = a_q = 1] &\gg \sum_{\substack{p,q \text{ distinct primes} \\ 5 \leq p,q \leq M}}   \frac{1}{p^{1/2 + \epsilon} q^{1/2 + \epsilon} M^2}  \\
		&= \frac{1}{M^2} \sum_{\substack{q \text{ prime} \\ 5 \leq q \leq M}} \frac{1}{ q^{1/2+\epsilon}} \sum_{\substack{p \text{ prime} \\ 5 \leq p \leq q}} \frac{1}{p^{1/2+\epsilon}}.
		\end{aligned}
	\end{equation} 
For all  $\epsilon_1, \epsilon_2 > 0$, we have 
	\begin{equation}\label{EqEst1}
	\begin{aligned}
	\sum_{\substack{p \text{ prime} \\ 5 \leq p \leq q}} \frac{1}{p^{1/2+\epsilon}} 
	&\sim \sum_{k = 2}^{\frac{q}{\log q}} \frac{1}{(k\log k)^{1/2+\epsilon}} \\
	&\gg \sum_{k = 2}^{\frac{q}{\log q}} \frac{1}{k^{1/2+\epsilon + \epsilon_1}} \\
	&\gg \int_{x=2}^{\frac{q}{\log q}} \frac{1}{x^{1/2+\epsilon+\epsilon_1}} dx \\
	&\gg q^{1/2-\epsilon-\epsilon_1-\epsilon_2}.
	\end{aligned}
	\end{equation}
Combining \eqref{EqProb1} and \eqref{EqEst1} yields
	\begin{align*}
		\Pr[a_p = a_q = 1] \gg \frac{1}{M^2} \sum_{\substack{q \text{ prime} \\5 \leq q \leq M}} \frac{q^{1/2-\epsilon-\epsilon_1-\epsilon_2} }{q^{1/2+\epsilon}} = \frac{1}{M^2} \sum_{\substack{q \text{ prime} \\5 \leq q \leq M}} \frac{ 1 }{ q^{2\epsilon + \epsilon_1 + \epsilon_2}}.
	\end{align*}
By replacing $2\epsilon + \epsilon_1 + \epsilon_2$ with $\epsilon$, we have
	\begin{align*}
	\Pr[a_p = a_q = 1] \gg \frac{1}{M^2} \sum_{\substack{q \text{ prime} \\5 \leq q \leq M}} \frac{ 1 }{ q^{\epsilon}}
	\end{align*}
	for all $\epsilon > 0$. Proceeding as in \eqref{EqEst1}, we bound $\Pr[a_p = a_q = 1]$ as
	\begin{align*}
		\Pr[a_p = a_q = 1] \gg \frac{1}{M^{1+\epsilon}}
	\end{align*}
	for all $\epsilon > 0$.
\end{proof}
\subsection{The proportion of choices for $p, q, E$ such that $p$ and $q$ are not anomalous primes for $E$} \label{SectionHard}
In this section, we find an upper bound for the probability 
$$\Pr[a_p \text{ or } a_q \neq 1 \text{ and } (p+1-a_p), (q+1-a_q) \mid (N+1-a_N) ].$$
Lemma \ref{LemmaUBound} identifies the upper bound by splitting into several cases. We can express the probability as a sum in which each summand corresponds to these cases and use Lemmas \ref{Lemmaapaq1} through \ref{LemmaFixPAP} to bound the summands.
\begin{lem} \label{Lemmaapaq1}
Let $5 \leq p < q$ be primes and assume that  $|a_p| \leq 2\sqrt{p}$, $|a_q| \leq 2\sqrt{q}$ and $(q+1-a_q) \mid (pq+1-a_pa_q)$. Then $a_q\neq 0$. Moreover, if not both of $a_p$ and $a_q$ are equal to $1$, then $a_q \neq1$.
\begin{proof}
Suppose for contradiction that $a_q = 0$. In particular, $(q+1) \mid (pq+1)$. Moreover, $q+1$ divides $pq+p$, so $q+1$ must divide $(pq+p) - (pq+1) = p-1$, but $0 < p-1 < q+1$. Hence, $a_q\neq 0$. Suppose for contradiction that $a_q = 1$. Here, $q \mid (1-a_p)$, but  
	\begin{align*}
	|1-a_p| \leq 1 + |a_p| \leq 1 + 2\sqrt{p} \leq 1 + 2\sqrt{q}.
	\end{align*}
Since  $q \geq 7$ we have $q>1+2\sqrt{q}$. Therefore $1-a_p = 0$, which contradicts that $a_p$ and $a_q$ are not both $1$. Hence, $a_q\neq 1$. 
\end{proof}
\end{lem}	
\begin{lem} \label{LemmaDiv}
Let $p$, $a_p$, $q$, and $a_q$ be integers. The divisibility conditions $(p+1-a_p), (q+1-a_q) \mid (pq+1-a_pa_q)$ hold if and only if 
	\begin{align*}
		(p+1-a_p) \mid (1-a_pa_q - q + qa_p) \qquad \text{and} \qquad (q+1-a_q) \mid (1-a_pa_q - p + pa_q).
	\end{align*}
	\begin{proof}
Suppose that $(p+1-a_p)$ divides $pq+1-a_pa_q$. This implies that $(p+1-a_p) \mid (1-a_pa_q - q + qa_p)$. It can be shown that $(p+1-a_p) \mid (1-a_pa_q - q + qa_p)$ implies $(p+1-a_p) \mid (pq+1-a_pa_q)$. Similarly, $(q+1-a_q) \mid (1-a_pa_q - p + pa_q)$ if and only if $(q+1-a_q) \mid (1-a_pa_q - p + pa_q)$.
	\end{proof}
\end{lem}

Considering Lemma \ref{LemmaDiv}, we will now refer to the divisibility conditions
$$ (p+1-a_p), (q+1-a_q) \mid (pq+1-a_pa_q)$$
interchangeably with
$$(p+1-a_p) \mid (1-a_pa_q - q + qa_p) \qquad \text{and} \qquad (q+1-a_q) \mid (1-a_pa_q - p + pa_q).$$

\begin{lem} \label{LemmaSolns}
Suppose that $p_0$ and $a_{p_0}$ are integers such that $(q+1-a_q) \mid (1 - a_{p_0} a_q - p_0 + p_0 a_q)$. If $p$ and $a_p$ are integers such that $(q+1-a_q) \mid (1-a_p a_q - p + p a_q)$, then $a_p = a_{p_0} + k(q+1-a_q) + (1-a_q) \alpha$  and $p = p_0 + k(q+1-a_q) - a_q \alpha$ for some integers $k$ and $\alpha$. Moreover, 
	\begin{align*}
		(1 - a_{p_0} a_q - p_0 + p_0 a_q) - (1-a_p a_q - p + p a_q) = k(q+1-a_q). 
	\end{align*}
	\begin{proof}
		Since $q+1-a_q$ divides both, $1 - a_{p_0} a_q - p_0 + p_0 a_q$ and $1-a_p a_q - p + p a_q$, it must divide 
		\begin{align*}
			(1 - a_{p_0} a_q - p_0 + p_0 a_q) - (1-a_p a_q - p + p a_q) = a_q (a_p - a_{p_0}) + (1-a_q) (p-p_0),
		\end{align*}
Let $x = a_p - a_{p_0}$ and $y = p-p_0$, so that
		\begin{align*}
			k(q+1-a_q) = a_q x + (1-a_q) y.
		\end{align*}
With $k$ fixed, this is a linear Diophantine equation in two variables. Since $\gcd(a_q, 1-a_q)=1$, all of the solutions take the form 
		\begin{align*}
			x &= k(q+1-a_q) + (1-a_q) \alpha \text{ and}\\
			y &= k(q+1-a_q) - a_q \alpha
		\end{align*}
		where $\alpha$ is an integer. 
	\end{proof}
\end{lem}

\begin{lem} \label{LemmaBoundPair1}
Let $q>6$ be a prime and  $a_q \neq 0, 1$ be an integer satisfying $\vert a_q\vert \leq 2\sqrt{q}$. For a prime $5 \leq p < q$ and $a_p$ an integer with $\vert a_p\vert \leq 2\sqrt{p}$, the number of distinct integer values of 
$\frac{ 1 - a_p a_q - p + p a_q}{q+1-a_q}$ is $O(\vert a_q\vert)$.
\begin{proof}
Let $p_0$ be a prime such that $5 \leq p_0 < q$ and $a_{p_0}$ be an integer such that $|a_{p_0}| \leq 2\sqrt{p_0}$ and $(q+1-a_q) \mid (1-a_{p_0} a_q - p_0 + p_0 a_q)$. Suppose that $p$ is a prime such that $5 \leq p < q$ and $a_{p}$ is an integer such that $|a_p| \leq 2\sqrt{p}$ and $(q+1-a_q) \mid (1-a_p a_q - p + pa_q)$. 
By Lemma \ref{LemmaSolns}, there exist integers $k$ and $\alpha$ such that $a_p = a_{p_0} + k(q+1-a_q) + (1-a_q) \alpha$ and $p = p_0 + k(q+1-a_q) - a_q \alpha$.
Compute $(1-a_{p_0} a_q - p_0 + p_0 a_q) - (1-a_p a_q - p + pa_q) = k(q+1-a_q)$. Thus, each value of $k$ corresponds to its integer value of $\frac{1-a_pa_q -p + p a_q}{q+1-a_q}$. Suppose that $\vert k \vert \geq 12 \vert a_q\vert$. Since $p<q$  

				$0 < p_0 + k(q+1-a_q) - a_q \alpha < q$,  and so 

	\begin{equation} \label{EqBound1}
	\begin{aligned}
		-p_0 - k(q+1-a_q) < - a_q \alpha < -p_0 - k(q+1-a_q) + q.
	\end{aligned}
	\end{equation}
	Adding $a_{p_0} + k(q+1-a_q) + \alpha$ to the above inequality, we have
	\begin{align*}
		a_{p_0} + \alpha - p_0 < a_{p_0} + k(q+1-a_q) + \alpha - a_q \alpha < a_{p_0} + \alpha - p_0 + q.
	\end{align*}
	Thus,
	\begin{equation} \label{EqBound2}
	\begin{aligned}
		a_{p_0} + \alpha - p_0 < a_p < a_{p_0} + \alpha - p_0 + q.
	\end{aligned}
	\end{equation}
Note that since $q \geq 7$ we have $3(q+1-a_q) > q$ and since $\vert k \vert > 12 \vert a_q\vert$ we have 
	\begin{align*}
		\left| \frac{k(q+1-a_q)}{a_q} \right| > 12(q+1-a_q) > 4q.
	\end{align*}
Moreover, since $0 < p_0 < q$, the quantities $\left\vert \frac{p_0}{a_q} \right\vert$ and $\left| \frac{p_0-q}{a_q} \right|$ are both at most $q$. In the case when $a_q > 0$, (\ref{EqBound1}) yields
	\begin{align*}
		\frac{p_0}{a_q} + \frac{k(q+1-a_q)}{a_q}  > \alpha > \frac{p_0-q}{a_q} + \frac{k(q+1-a_q)}{a_q}.
	\end{align*}
If $k > 0$ as well, then $\frac{k(q+1-a_q)}{a_q} > 0$, and so
	\begin{align*}
	\alpha > \frac{p_0-q}{a_q} + \frac{k(q+1-a_q)}{a_q} >  3q.
	\end{align*}
Since $|a_{p_0}| < 2\sqrt{p_0} < 2\sqrt{q} < q$, (\ref{EqBound2}) implies that
	\begin{align*}
		q = -q+3q-q < a_{p_0} + \alpha - p_0 < a_p,
	\end{align*}
which is a contradiction.  If $k < 0$, then $\frac{k(q+1-a_q)}{a_q} < 0$, and so 
	\begin{align*}
	\alpha < \frac{p_0}{q} + \frac{k(q+1-a_q)}{a_q} <  -3q.
	\end{align*}
This time, (\ref{EqBound2}) yields
	\begin{align*}
	a_p < a_{p_0} + \alpha - p_0 + q < q - 3q - 0 + q = -q,
	\end{align*}
but this is a contradiction. Now assume that $a_q < 0$. By (\ref{EqBound1}), 
	\begin{align*}
		\frac{p_0}{a_q} + \frac{k(q+1-a_q)}{a_q} < \alpha <  \frac{p_0-q}{a_q} + \frac{k(q+1-a_q)}{a_q}.
	\end{align*}
If $k > 0$, then $\frac{k(q+1-a_q)}{a_q} < 0$, and so 
	\begin{align*}
		\alpha < \frac{p_0-q}{a_q} + \frac{k(q+1-a_q)}{a_q} <  -3q.
	\end{align*}
Therefore, (\ref{EqBound2}) implies $a_p < a_{p_0} + \alpha - p_0 + q < -q$. If $k < 0$, then $\frac{k(q+1-a_q)}{a_q} > 0$, and so
	\begin{align*}
		\alpha > \frac{p_0}{a_q} + \frac{k(q+1-a_q)}{a_q} >  3q.
	\end{align*}
Again, (\ref{EqBound2}) implies that $a_p > a_{p_0} + \alpha - p_0 > q$. In all cases, $\vert a_p\vert > q$ as desired. Hence, $\vert k\vert \leq 12\vert a_q\vert$, which means that the number of possible distinct values of $k$ and, by extension, the number of possible distinct integer values of $\frac{1-a_pa_q-p+pa_q}{q+1-a_q}$ is $O(a_q)$. 
	\end{proof}
\end{lem}

\begin{lem}\label{LemmaDivBound}
Let $n$ be a positive integer. The number of divisors $d(n)$ of $n$ satisfies $d(n) = o(n^\epsilon)$ for all $\epsilon > 0$. 
\begin{proof}
See \cite[page 296, Theorem 13.12]{Apostol}.
\end{proof}
\end{lem}

\begin{lem} \label{LemmaBoundPair2}
Let $q\geq 17 $ be a prime and let $a_q \neq 0,1$ be an integer satisfying $9 < \vert a_q\vert \leq 2\sqrt{q}$. Let $p$ be a prime such that $5 \leq p < q$ and let $a_p$ be an integer satisfying $\vert a_p\vert \leq 2\sqrt{p}$. For a fixed integer
	\begin{align*}
		l_0= \frac{1-a_pa_q - p + pa_q}{q+1-a_q},
	\end{align*}
the number of distinct pairs $(p,a_p)$ with $(p+1-a_p) \mid (1-a_pa_q - q + qa_p)$ is $o(q^\epsilon)$ for all $\epsilon > 0$. 
	\begin{proof}
Suppose that $p_0$ is a prime such that $5 \leq p_0 < q$ and $a_{p_0}$ is an integer satisfying $|a_{p_0}| \leq 2 \sqrt{p_0}$ such that
		\begin{align*}
			l_0 = \frac{1-a_{p_0} a_q - p_0 + p_0 a_q}{q+1-a_q}
		\end{align*}
and $(p_0+1-a_{p_0}) \mid (1-a_{p_0} a_q - q + qa_{p_0})$. Further suppose that $p$ is a prime such that $5 \leq p < q$ and  $|a_p| \leq 2\sqrt{p}$ is an integer such that
		\begin{align*}
			l_0 = \frac{1-a_p a_q - p + p a_q}{q+1-a_q}
		\end{align*}
and $(p+1-a_p) \mid (1-a_p a_q - q a_p)$. By Lemma \ref{LemmaSolns}, there are integers $k$ and $\alpha$ such that $a_p = a_{p_0} + k(q+1-a_q) + (1-a_q) \alpha$ and $p = p_0 + k(q+1-a_q) -a_q \alpha$.
Note that $(1-a_{p_0} a_q - p_0 + p_0 a_q) - (1-a_p a_q - p + pa_q) = k(q+1-a_q)$. Since 
		\begin{align*}
			\frac{1-a_{p_0} a_q - p_0 + p_0 a_q}{q+1-a_q} =  \frac{1-a_p a_q - p + p a_q}{q+1-a_q}
		\end{align*}
		
\flushleft we have that $k = 0$. In particular, $\alpha$ is $O(q)$ because $0 < p < q$. Compute
		\begin{align*}
			p+1-a_p &= (p_0-a_q \alpha) + 1 - (a_{p_0} + (1-a_q) \alpha) \\
						&= p_0 + 1 - a_{p_0} - \alpha
		\end{align*}
and
		\begin{align*}
			1-a_pa_q-q+qa_p &= 1-a_pa_q - q(1-a_p) \\											
			&= 1 - a_{p_0} a_q - q + qa_{p_0} + (q-a_q)(1-a_q) \alpha.
		\end{align*}
Let $d = p_0+1-a_{p_0}$ and let $n = 1-a_{p_0}a_q - q + qa_{p_0}$ such that $d \mid n$. Moreover, $p+1-a_p = d - \alpha$, $1-a_pa_q - q + qa_p = n + (q-a_q)(1-a_q) \alpha$, and $(d-\alpha) \mid (n+(q-a_q)(1-a_q) \alpha)$. Note that
		\begin{align*}
			\frac{n}{d} - \frac{n+(q-a_q)(1-a_q) \alpha}{d-\alpha}
		\end{align*}
is an integer. Compute
		\begin{align*}
			\frac{n}{d} - \frac{n+(q-a_q)(1-a_q)\alpha}{d-\alpha} &= \frac{n(d-\alpha) - d(n+(q-a_q)(1-a_q)\alpha)}{d(d-\alpha)} \\
												      &= \frac{ -\frac{n}{d} \alpha - (q-a_q)(1-a_q)\alpha}{d-\alpha},
		\end{align*}
so $(d-\alpha) \mid \left( -\frac{n}{d} - (q-a_q)(1-a_q) \right) \alpha$. Thus,
		\begin{align*}
			\frac{d - \alpha}{\gcd(d-\alpha,\alpha)} \mid \left( -\frac{n}{d} - (q-a_q)(1-a_q) \right).
		\end{align*}
Since $\gcd(d-\alpha,\alpha) = \gcd(d,\alpha)$,  $\frac{d - \alpha}{\gcd(d,\alpha)} \mid \left( -\frac{n}{d} - (q-a_q)(1-a_q) \right)$.
	
Whenever $\alpha$ satisfies the above divisibility condition, there is some $d'$ dividing $-\frac{n}{d} - (q-a_q)(1-a_q)$ such that $\frac{d-\alpha}{\gcd(d,\alpha)} = d'$, or equivalently $d-\alpha = d' \gcd(d,\alpha)$.
	
Similarly there is $g \mid d$ such that $\alpha = d-d'g$. Since $d = p_0 + 1 - a_{p_0}$ and $5 \geq p_0$, we have that $d\neq 0$. We will show that $-\frac{n}{d} - (q-a_q)(1-a_q)\neq 0$. Note that $p_0 < 3(p_0+1-a_{p_0})$ and $1+|a_{p_0}| \leq p_0$.  Then, $\left| \frac{n}{d}\right| < \frac{3(1+|a_{p_0} a_q| + p_0q)}{p_0}$. Moreover, $1+|a_{p_0} a_q| \leq 1+4\sqrt{p_0q} < \frac{p_0q}{2}$ because $p_0q \geq 5 \cdot 17 = 85$, and so
		\begin{align*}
			\left| \frac{n}{d}\right| &< \frac{3 \left( \frac{p_0q}{2} + p_0q \right)}{p_0} = \frac{9}{2}q.
		\end{align*}
On the other hand, we have $q^2 > 16q \geq 4a_q^2$ which implies that $q-a_q > \frac{q}{2}$. Since $|a_q| > 9$, 		
		\begin{align*}
			\left| \frac{n}{d} \right| &< \frac{9}{2} q \leq \frac{1}{2} q |1-a_q| < (q-a_q)|1-a_q| = |(q-a_q)(1-a_q)|.
		\end{align*}
Therefore, $-\frac{n}{d} - (q-a_q) (1-a_q)\neq 0$ as desired. Note that $n,d$, and $-\frac{n}{d}-(q-a_q)(1-a_q)$ are all fixed with respect to $q,a_q, p_0,$ and $a_{p_0}$. Their bounds are $n = o(q\sqrt{p})$, $d = o(p)$ and $-\frac{n}{d} - (q-a_q)(1-a_q) = o(q\sqrt{p})$. Therefore, $d = o(q)$ and $-\frac{n}{d} - (q-a_q)(1-a_q) = o(q^2)$. By Lemma \ref{LemmaDivBound}, there are $o(q^\epsilon)$ and $o(q^{2\epsilon})$ possible values of $d'$ and $g$ for all $\epsilon > 0$ respectively, so there are $o(q^{3\epsilon})$ possible values of $\alpha$. Consequently, there are  $o(q^{\epsilon})$ possible combinations of $(p,a_p)$ for all $\epsilon > 0$.
\end{proof}
\end{lem}

\begin{lem} \label{LemmaFixQAQ}
Let $q>6$ be a prime $a_q \neq 1$ be an integer satisfying $|a_q| \leq 2\sqrt{q}$. The following hold
\begin{enumerate}
	\item For a fixed integer $a_p$ there are $O(1)$ integers $p$ with $5 \leq p < q$ satisfying 
		\begin{align*}
		(q+1-a_q) \mid (1-a_pa_q - p + pa_q).
	\end{align*}	
	\item For a fixed integer $p$  with $0 < p < q$ and given $a_q = O(1)$, there are $O(1)$ integers $a_p$ with $|a_p| \leq 2\sqrt{p}$ satisfying
	\begin{align*}
		(q+1-a_q) \mid (1-a_pa_q - p + pa_q).
	\end{align*}
\end{enumerate}

\begin{proof}
\emph{(1)} Note that $q+1-a_q$ and $1-a_pa_q$ are fixed. Furthermore, 
	\begin{align*}
		1-a_pa_q - p + pa_q = 1-a_pa_q - p(1-a_q).
	\end{align*}
Let $p_0$ and $p$ be two integers with $5 \leq p,p_0 < q$ satisfying
	\begin{align*}
		(q+1-a_q) \mid (1-a_pa_q - p+pa_q), (1-a_pa_q - p_0 + p_0 a_q).
	\end{align*}
In particular,
	\begin{align*}
			0 &\equiv (1-a_pa_q - p+pa_q) - (1-a_pa_q - p_0 + p_0 a_q) \\
			   &\equiv (p_0-p)(1-a_q) \mod{(q+1-a_q)},
	\end{align*}
or equivalently, $(q+1-a_q) \mid (p_0-p)(1-a_q)$. Since $a_q \neq 1$, $\gcd(q+1-a_q, 1-a_q) = \gcd(q,1-a_q) = 1$. Therefore, $(q+1-a_q) \mid (p_0-p)$. However, $q+1-a_q = \Theta(q)$, but $5 \leq p,p_0 < q$, so there are $O(1)$ possible values of $p$ satisfying
	\begin{align*}
		(q+1-a_q) \mid (1-a_pa_q - p + pa_q).
	\end{align*}
	
\emph{(2)} Note that $q+1-a_q$ and $1 - p + pa_q$ are fixed. Suppose that $a_p$ and $a_{p_0}$ are both integers with $|a_p|,|a_{p_0}| \leq 2\sqrt{p}$ and 
	\begin{align*}
		(q+1-a_q) \mid (1-a_pa_q  - p + pa_q), (1-a_{p_0} a_q - p + pa_q).
	\end{align*}
In particular,
	\begin{align*}
		        0 
			&\equiv (1-a_pa_q  - p + pa_q) - (1-a_{p_0} a_q - p + pa_q) \\
			&\equiv (a_{p_0} - a_p)a_q \mod {(q+1-a_q)},
	\end{align*}
Thus, $\frac{q+1-a_q}{\gcd(q+1-a_q,a_q)}$ divides $(a_{p_0}-a_p)$.	 Since $a_q$ is $O(1)$, so is $\gcd(q+1-a_q,a_q)$ and thus $\frac{q+1-a_q}{\gcd(q+1-a_q,a_q)}$ is $\Theta(q)$. However, $|a_p|, |a_{p_0}| \leq 2\sqrt{p} < 2\sqrt{q}$, so there are $O(1)$ possible values of $a_p$ as desired. 	

\end{proof}
\end{lem}
\begin{lem}\label{LemmaFixPAP}
Let $5 \leq p \leq 13$ be a prime and $a_p$ be an integer satisfying $|a_p| \leq 2\sqrt{p}$.
	\begin{enumerate}
		\item  There are $O(1)$ possible values of $a_q$ satisfying $1-a_pa_q - p + pa_q = 0$.
		\item  For an integer $a_q$ with $1-a_pa_q - p + pa_q \neq 0$ there are $O(1)$ integers $q$ satisfying $|a_q| \leq 2\sqrt{q}$ and $(q+1-a_q) \mid (1-a_pa_q - p + pa_q)$.
	\end{enumerate}
	\begin{proof}
{(1)} Note that $a_p-p \neq 0$ because $p \geq 5$. If $1-a_pa_q - p + pa_q = 0$, then $\frac{1-p}{a_p-p} = a_q$. Since $p = O(1)$, $a_p = O(1)$ as well. Thus, there are $O(1)$ possible values of $a_q$. 

{(2)} Again, since $p = O(1)$, $a_p = O(1)$ as well. Therefore, $1-a_pa_q - p + pa_q = O(a_q)$, but $q+1-a_q = \Theta(q)$. Thus there are $O(1)$ possible values of $q$ satisfying $(q+1-a_q) \mid (1-a_pa_q - p + pa_q)$ and by extension, $O(1)$ integers $q$ with $|a_q| \leq 2\sqrt{q}$ satisfying the divisibility condition.

\end{proof}
\end{lem}
\begin{lem}\cite[Corollary 4.8]{REU2016} \label{LemmaAnomProd}
Let $E/\Q$ be an elliptic curve and let $N = pq$ be an elliptic Korselt number of Type I for $E$ such that $p < q$. One of the following holds:
	\begin{enumerate}
		\item  $p \leq 13$
		\item  $p$ and $q$ are anomalous for $E$.
		\item  $p \geq \frac{\sqrt{q}}{16}$. 
	\end{enumerate}	
\end{lem}
\begin{lem} \label{LemmaUBound}
Let $N = pq$ be a product of two distinct primes $5 \leq p,q \leq M$ chosen at random. Let $\Emod{N}$ be an elliptic curve with good reduction at $p$ and $q$. The probability that $(p+1-a_p), (q+1-a_q) \mid (N+1-a_N)$ and $a_p$ and $a_q$ are not both $1$ is 
	$$ O\left( \frac{1}{M^{5/4-\epsilon}} \right)$$
for all $\epsilon > 0$. 
\begin{proof}
\ifshowold{	
Let $a_p$ and $a_q$ be integers such that $|a_p| \leq 2\sqrt{p}$ and $|a_q| \leq 2\sqrt{q}$ and let $a_N = a_pa_q$. Note that the number of pairs of $p$ and $q$ is on the order of  $\left( \frac{M}{\log M} \right)^2$. Also, note that
	\begin{equation} \label{EqIneq14}
	\begin{aligned}
		 &\Pr \left[\substack{(p+1-a_p),(q+1-a_q) \mid (N+1-a_N), \\
					a_p \text{ or } a_q \neq 1} \right] \\
		&= \sum_{\substack{ q,a_q,p,a_p \\ (p+1-a_p),(q+1-a_q) \mid (N+1-a_N)  \\ a_p \text{ or } a_q \neq 1} } \Pr\left[ \substack{ p,q \text{ chosen}, \\ \#\Emod{p} = p +1 - a_p, \\ \#\Emod{q} = q+1-a_q} \right] \\
		&= \sum_{\substack{ q,a_q,p,a_p \\ (p+1-a_p),(q+1-a_q) \mid (N+1-a_N)  \\ a_p \text{ or } a_q \neq 1} } \Pr\left[ \substack{ \#\Emod{p} = p +1 - a_p, \\ \#\Emod{q} = q+1-a_q} \biggr\vert p,q \text{ chosen} \right]  \\
		&\phantom{=} \cdot \left( \frac{1}{\#\{(p,q) \mid p,q \text{ distinct primes}, 5 \leq p,q \leq M\}} \right) \\
		&\approx \left( \frac{\log M}{M} \right)^2 \sum_{\substack{q,a_q,p,a_p \\ (p+1-a_p),(q+1-a_q) \mid (N+1-a_N)  \\ a_p \text{ or } a_q \neq 1} } \Pr\left[ \substack{ \#\Emod{p} = p +1 - a_p, \\ \#\Emod{q} = q+1-a_q} \biggr\vert p,q \text{ chosen} \right].
	\end{aligned}
	\end{equation}
Furthermore,
	\begin{align*}
		 &\sum_{\substack{q,a_q,p,a_p \\ (p+1-a_p),(q+1-a_q) \mid (N+1-a_N)  \\ a_p \text{ or } a_q \neq 1} } \Pr\left[ \substack{ \#\Emod{p} = p +1 - a_p, \\ \#\Emod{q} = q+1-a_q} \biggr\vert p,q \text{ chosen} \right] \\
		&=2 \cdot \left( \sum_{\substack{q,a_q,p,a_p \\ (p+1-a_p),(q+1-a_q) \mid (N+1-a_N)  \\ a_p \text{ or } a_q \neq 1 \\ p < q} } \Pr\left[ \substack{ \#\Emod{p} = p +1 - a_p, \\ \#\Emod{q} = q+1-a_q} \biggr\vert p,q \text{ chosen} \right] \right) \\
		&\approx \sum_{\substack{q,a_q,p,a_p \\ (p+1-a_p),(q+1-a_q) \mid (N+1-a_N)  \\ a_p \text{ or } a_q \neq 1 \\ p < q} } \Pr\left[ \substack{ \#\Emod{p} = p +1 - a_p, \\ \#\Emod{q} = q+1-a_q} \biggr\vert p,q \text{ chosen} \right].
\end{align*}
By Lemma \ref{Lemmaapaq1},

\begin{align*}
	&\sum_{\substack{ q,a_q,p,a_p \\ (p+1-a_p),(q+1-a_q) \mid (N+1-a_N)  \\ a_p \text{ or } a_q \neq 1 \\ p < q} } \Pr\left[ \substack{ \#\Emod{p} = p +1 - a_p, \\ \#\Emod{q} = q+1-a_q} \biggr\vert p,q \text{ chosen} \right] \\
	&= \sum_{\substack{ q,a_q,p,a_p \\ (p+1-a_p),(q+1-a_q) \mid (N+1-a_N)  \\ a_q \neq 1 \\ p < q} } \Pr\left[ \substack{ \#\Emod{p} = p +1 - a_p, \\ \#\Emod{q} = q+1-a_q} \biggr\vert p,q \text{ chosen} \right].
\end{align*}
Moreover,

\begin{equation}\label{EqIneq12}
\begin{aligned}
		&\sum_{\substack{ q,a_q,p,a_p \\ (p+1-a_p),(q+1-a_q) \mid (N+1-a_N)  \\ a_q \neq 1 \\ p < q} } \Pr\left[ \substack{ \#\Emod{p} = p +1 - a_p, \\ \#\Emod{q} = q+1-a_q} \biggr\vert p,q \text{ chosen} \right] \\
		&= \sum_{\substack{q,a_q \\ q \geq 17 \\ 9 < |a_q| \leq 2\sqrt{q}}} \sum_{\substack{p,a_p \\ p < q \\ (p+1-a_p) \mid (N+1-a_N) \\ (q+1-a_q) \mid (N+1-a_N)}} \Pr\left[ \substack{ \#\Emod{p} = p +1 - a_p, \\ \#\Emod{q} = q+1-a_q} \biggr\vert p,q \text{ chosen} \right] \\
		&+ \sum_{q \geq 17} \sum_{\substack{ |a_q| \leq 9  \\ a_q \neq 1} } \sum_{\substack{p,a_p \\ p < q \\ (p+1-a_p) \mid (N+1-a_N) \\ (q+1-a_q) \mid (N+1-a_N)}} \Pr\left[ \substack{ \#\Emod{p} = p +1 - a_p, \\ \#\Emod{q} = q+1-a_q} \biggr\vert p,q \text{ chosen} \right] \\
		&+ \sum_{\substack{q,a_q,p,a_p \\ (p+1-a_p),(q+1-a_q) \mid (N+1-a_N)  \\ a_q \neq 1 \\ p < q < 17} } \Pr\left[ \substack{ \#\Emod{p} = p +1 - a_p, \\ \#\Emod{q} = q+1-a_q} \biggr\vert p,q \text{ chosen} \right].
	\end{aligned}
\end{equation}
If $p < q < 17$, then there are only $O(1)$ possible combinations of $q,a_q,p$ and $a_p$. Thus,
	\begin{equation}\label{EqIneq13}
	\begin{aligned}
		\sum_{\substack{q,a_q,p,a_p \\ (p+1-a_p),(q+1-a_q) \mid (N+1-a_N)  \\ a_p \text{ or } a_q \neq 1 \\ p < q < 17} } \Pr\left[ \substack{ \#\Emod{p} = p +1 - a_p, \\ \#\Emod{q} = q+1-a_q} \biggr\vert p,q \text{ chosen} \right] = O(1).
	\end{aligned}
	\end{equation}
By Lemma \ref{LemmaAnomProd}, 
	\begin{equation}\label{EqIneq1}
	\begin{aligned} 
		&\sum_{q \geq 17} \sum_{\substack{ |a_q| \leq 9  \\ a_q \neq 1} } \sum_{\substack{p,a_p \\ p < q \\ (p+1-a_p) \mid (N+1-a_N) \\ (q+1-a_q) \mid (N+1-a_N)}} \Pr\left[ \substack{ \#\Emod{p} = p +1 - a_p, \\ \#\Emod{q} = q+1-a_q} \biggr\vert p,q \text{ chosen} \right] \\
		&= \sum_{q \geq 17} \sum_{\substack{ |a_q| \leq 9  \\ a_q \neq 1} } \sum_{\substack{p,a_p \\ p \leq 13 \\ (p+1-a_p) \mid (N+1-a_N) \\ (q+1-a_q) \mid (N+1-a_N)}} \Pr\left[ \substack{ \#\Emod{p} = p +1 - a_p, \\ \#\Emod{q} = q+1-a_q} \biggr\vert p,q \text{ chosen} \right] \\
		&+ \sum_{q \geq 17} \sum_{\substack{ |a_q| \leq 9  \\ a_q \neq 1} } \sum_{\substack{p,a_p \\ 13 < p < q \\ (p+1-a_p) \mid (N+1-a_N) \\ (q+1-a_q) \mid (N+1-a_N)}} \Pr\left[ \substack{ \#\Emod{p} = p +1 - a_p, \\ \#\Emod{q} = q+1-a_q} \biggr\vert p,q \text{ chosen} \right].
	\end{aligned}
	\end{equation}
By Corollary \ref{corProb} the first summand of (\ref{EqIneq1}) satisfies
	\begin{align*}
		&\sum_{q \geq 17} \sum_{\substack{ |a_q| \leq 9  \\ a_q \neq 1} } \sum_{\substack{p,a_p \\ p \leq 13 \\ (p+1-a_p) \mid (N+1-a_N) \\ (q+1-a_q) \mid (N+1-a_N)}} \Pr\left[ \substack{ \#\Emod{p} = p +1 - a_p, \\ \#\Emod{q} = q+1-a_q} \biggr\vert p,q \text{ chosen} \right] \\
		&\ll \sum_{q \geq 17} \sum_{\substack{ |a_q| \leq 9  \\ a_q \neq 1} } \sum_{\substack{p,a_p \\ p \leq 13 \\ (p+1-a_p) \mid (N+1-a_N) \\ (q+1-a_q) \mid (N+1-a_N)}} \frac{1}{(pq)^{1/2-\epsilon}} \\
		&\ll \sum_{q \geq 17} \sum_{\substack{ |a_q| \leq 9  \\ a_q \neq 1} } \sum_{\substack{p,a_p \\ p \leq 13 \\ (p+1-a_p) \mid (N+1-a_N) \\ (q+1-a_q) \mid (N+1-a_N)}} \frac{1}{q^{1/2-\epsilon}}.
	\end{align*}
There are $O(1)$ possible combinations of $a_q,p$ and $a_p$ in the above sum, so
	\begin{equation} \label{EqIneq2}
	\begin{aligned}
		&\sum_{q \geq 17} \sum_{\substack{ |a_q| \leq 9  \\ a_q \neq 1} } \sum_{\substack{p,a_p \\ p \leq 13 \\ (p+1-a_p) \mid (N+1-a_N) \\ (q+1-a_q) \mid (N+1-a_N)}} \Pr\left[ \substack{ \#\Emod{p} = p +1 - a_p, \\ \#\Emod{q} = q+1-a_q} \biggr\vert p,q \text{ chosen} \right] \\
		&\ll \sum_{q \geq 17} \frac{1}{q^{1/2-\epsilon}} \\
		&\ll \int_{17}^{M} \frac{1}{x^{1/2-\epsilon}} dx \\
		&\ll M^{1/2+\epsilon}.
	\end{aligned}
	\end{equation}
On the other hand, $|a_p| \leq 2\sqrt{p} < 2\sqrt{q}$, so Lemma \ref{LemmaNumAp} and Lemma \ref{LemmaBounds} show that the second summand of (\ref{EqIneq1}) satisfies
	\begin{align*}
		 &\sum_{q \geq 17} \sum_{\substack{ |a_q| \leq 9  \\ a_q \neq 1} } \sum_{\substack{p,a_p \\ 13 < p < q \\ (p+1-a_p) \mid (N+1-a_N) \\ (q+1-a_q) \mid (N+1-a_N)}} \Pr\left[ \substack{ \#\Emod{p} = p +1 - a_p, \\ \#\Emod{q} = q+1-a_q} \biggr\vert p,q \text{ chosen} \right] \\
		&\ll  \sum_{q \geq 17} \sum_{\substack{ |a_q| \leq 9  \\ a_q \neq 1} } \sum_{|a_p| \leq 2\sqrt{q}} \sum_{\substack{p \\ 13 < p < q \\ (p+1-a_p) \mid (N+1-a_N) \\ (q+1-a_q) \mid (N+1-a_N)}} \Pr\left[ \substack{ \#\Emod{p} = p +1 - a_p, \\ \#\Emod{q} = q+1-a_q} \biggr\vert p,q \text{ chosen} \right] \\
		&\ll \sum_{q \geq 17} \sum_{\substack{ |a_q| \leq 9  \\ a_q \neq 1} } \sum_{|a_p| \leq 2\sqrt{q}} \sum_{\substack{p \\ 13 < p < q \\ (p+1-a_p) \mid (N+1-a_N) \\ (q+1-a_q) \mid (N+1-a_N)}} \frac{1}{(pq)^{1/2-\epsilon}}.
	\end{align*}
By lemma \ref{LemmaAnomProd}, $p \geq \frac{\sqrt{q}}{16}$ in the above sum. Thus,
	\begin{align*}
		 &\sum_{q \geq 17} \sum_{\substack{ |a_q| \leq 9  \\ a_q \neq 1} } \sum_{\substack{p,a_p \\ 13 < p < q \\ (p+1-a_p) \mid (N+1-a_N) \\ (q+1-a_q) \mid (N+1-a_N)}} \Pr\left[ \substack{ \#\Emod{p} = p +1 - a_p, \\ \#\Emod{q} = q+1-a_q} \biggr\vert p,q \text{ chosen} \right] \\
		&\ll \sum_{q \geq 17} \sum_{\substack{ |a_q| \leq 9  \\ a_q \neq 1} } \sum_{|a_p| \leq 2\sqrt{q}} \sum_{\substack{p \\ \frac{\sqrt{q}}{16} < p < q \\ (p+1-a_p) \mid (N+1-a_N) \\ (q+1-a_q) \mid (N+1-a_N)}} \frac{1}{(pq)^{1/2-\epsilon}} \\
		&\ll \sum_{q \geq 17} \sum_{\substack{ |a_q| \leq 9  \\ a_q \neq 1} } \sum_{|a_p| \leq 2\sqrt{q}} \sum_{\substack{p \\ \frac{\sqrt{q}}{16} < p < q \\ (p+1-a_p) \mid (N+1-a_N) \\ (q+1-a_q) \mid (N+1-a_N)}} \frac{1}{(q^{3/2})^{1/2-\epsilon}} \\
		&= \sum_{q \geq 17} \sum_{\substack{ |a_q| \leq 9  \\ a_q \neq 1} } \sum_{|a_p| \leq 2\sqrt{q}} \sum_{\substack{p \\ \frac{\sqrt{q}}{16} < p < q \\ (p+1-a_p) \mid (N+1-a_N) \\ (q+1-a_q) \mid (N+1-a_N)}} \frac{1}{q^{3/4-3\epsilon/2}}.
	\end{align*}
Lemma \ref{LemmaFixQAQ} implies that each combination of $q,a_q$ and $a_p$ yields only $O(1)$ possible values of $p$, so
	\begin{align*}
	&\sum_{q \geq 17} \sum_{\substack{ |a_q| \leq 9  \\ a_q \neq 1} } \sum_{\substack{p,a_p \\ 13 < p < q \\ (p+1-a_p) \mid (N+1-a_N) \\ (q+1-a_q) \mid (N+1-a_N)}} \Pr\left[ \substack{ \#\Emod{p} = p +1 - a_p, \\ \#\Emod{q} = q+1-a_q} \biggr\vert p,q \text{ chosen} \right] \\
	&\ll \sum_{q \geq 17} \sum_{\substack{ |a_q| \leq 9  \\ a_q \neq 1} } \sum_{|a_p| \leq 2\sqrt{q}} \frac{1}{q^{3/4-3\epsilon/2}}  \\
	&\ll \sum_{q \geq 17} \sum_{\substack{ |a_q| \leq 9  \\ a_q \neq 1} } \frac{1}{q^{1/4-3\epsilon/2}} \\
	&\ll \sum_{q \geq 17} \frac{1}{q^{1/4-3\epsilon/2}} \\
	&\ll \int_{17}^{M} \frac{1}{x^{1/4-3\epsilon/2}} dx \\
	&\ll M^{3/4+3\epsilon/2}.
	\end{align*}
Replacing $\epsilon$ with $2\epsilon/3$, we have that 
	\begin{equation} \label{EqIneq3}
	\begin{aligned}
		&\sum_{q \geq 17} \sum_{\substack{ |a_q| \leq 9  \\ a_q \neq 1} } \sum_{\substack{p,a_p \\ 13 < p < q \\ (p+1-a_p) \mid (N+1-a_N) \\ (q+1-a_q) \mid (N+1-a_N)}} \Pr\left[ \substack{ \#\Emod{p} = p +1 - a_p, \\ \#\Emod{q} = q+1-a_q} \biggr\vert p,q \text{ chosen} \right] \\\
		&\ll M^{3/4+\epsilon}
	\end{aligned}
	\end{equation}
for all $\epsilon > 0$. Combining (\ref{EqIneq1}), (\ref{EqIneq2}) and (\ref{EqIneq3}) shows that
	\begin{equation} \label{EqIneq4}
	\begin{aligned}
		&\sum_{q \geq 17} \sum_{\substack{ |a_q| \leq 9  \\ a_q \neq 1} } \sum_{\substack{p,a_p \\ p < q \\ (p+1-a_p) \mid (N+1-a_N) \\ (q+1-a_q) \mid (N+1-a_N)}} \Pr\left[ \substack{ \#\Emod{p} = p +1 - a_p, \\ \#\Emod{q} = q+1-a_q} \biggr\vert p,q \text{ chosen} \right] \\
		&\ll M^{3/4+\epsilon}.
	\end{aligned}
	\end{equation}
Using lemma \ref{LemmaAnomProd}, separate the first summand
	\begin{align*}
	 \sum_{\substack{q,a_q \\ q \geq 17 \\ 9 < |a_q| \leq 2\sqrt{q}}} \sum_{\substack{p,a_p \\ p < q \\ (p+1-a_p) \mid (N+1-a_N) \\ (q+1-a_q) \mid (N+1-a_N)}} \Pr\left[ \substack{ \#\Emod{p} = p +1 - a_p, \\ \#\Emod{q} = q+1-a_q} \biggr\vert p,q \text{ chosen} \right]
	\end{align*}
	of (\ref{EqIneq12}) into 
	\begin{equation} \label{EqIneq10}
	\begin{aligned}
		&\sum_{\substack{q,a_q \\ q \geq 17 \\ 9 < |a_q| \leq 2\sqrt{q}}} \sum_{\substack{p,a_p \\ p < q \\ (p+1-a_p) \mid (N+1-a_N) \\ (q+1-a_q) \mid (N+1-a_N)}} \Pr\left[ \substack{ \#\Emod{p} = p +1 - a_p, \\ \#\Emod{q} = q+1-a_q} \biggr\vert p,q \text{ chosen} \right]\\ 
		&= \sum_{\substack{q,a_q \\ q \geq 17 \\ 9 < |a_q| \leq 2\sqrt{q}}} \sum_{\substack{p,a_p \\ p \leq 13 \\ (p+1-a_p) \mid (N+1-a_N) \\ (q+1-a_q) \mid (N+1-a_N)}} \Pr\left[ \substack{ \#\Emod{p} = p +1 - a_p, \\ \#\Emod{q} = q+1-a_q} \biggr\vert p,q \text{ chosen} \right] \\
		&+ \sum_{\substack{q,a_q \\ q \geq 17 \\ 9 < |a_q| \leq 2\sqrt{q}}} \sum_{\substack{p,a_p \\ \frac{\sqrt{q}}{16} \leq p < q \\ (p+1-a_p) \mid (N+1-a_N) \\ (q+1-a_q) \mid (N+1-a_N)}} \Pr\left[ \substack{ \#\Emod{p} = p +1 - a_p, \\ \#\Emod{q} = q+1-a_q} \biggr\vert p,q \text{ chosen} \right].
	\end{aligned}
	\end{equation}	
Rearrange the first of the two summands in (\ref{EqIneq10}) to obtain
	\begin{align*}
		&\sum_{\substack{q,a_q \\ q \geq 17 \\ 9 < |a_q| \leq 2\sqrt{q}}} \sum_{\substack{p,a_p \\ p \leq 13 \\ (p+1-a_p) \mid (N+1-a_N) \\ (q+1-a_q) \mid (N+1-a_N)}} \Pr\left[ \substack{ \#\Emod{p} = p +1 - a_p, \\ \#\Emod{q} = q+1-a_q} \biggr\vert p,q \text{ chosen} \right] \\
		&=\sum_{\substack{p,a_p \\ p \leq 13}} \sum_{\substack{q,a_q \\ q \geq 17 \\ 9 < |a_q| \leq 2\sqrt{q} \\ (p+1-a_p) \mid (N+1-a_N) \\ (q+1-a_q) \mid (N+1-a_N)}}  \Pr\left[ \substack{ \#\Emod{p} = p +1 - a_p, \\ \#\Emod{q} = q+1-a_q} \biggr\vert p,q \text{ chosen} \right] \\
		&= \sum_{\substack{p,a_p \\ p \leq 13}} \sum_{\substack{q \geq 17}} \sum_{\substack{9 < |a_q| \leq 2\sqrt{q} \\ (p+1-a_p) \mid (N+1-a_N) \\ (q+1-a_q) \mid (N+1-a_N)}}  \Pr\left[ \substack{ \#\Emod{p} = p +1 - a_p, \\ \#\Emod{q} = q+1-a_q} \biggr\vert p,q \text{ chosen} \right] \\
		&= \sum_{\substack{p,a_p \\ p \leq 13}} \sum_{\substack{q \geq 17}} \sum_{\substack{9 < |a_q| \leq 2\sqrt{q} \\ (p+1-a_p) \mid (N+1-a_N) \\ (q+1-a_q) \mid (N+1-a_N) \\ 1-a_pa_q-p+pa_q = 0}}  \Pr\left[ \substack{ \#\Emod{p} = p +1 - a_p, \\ \#\Emod{q} = q+1-a_q} \biggr\vert p,q \text{ chosen} \right] \\
		&+ \sum_{\substack{p,a_p \\ p \leq 13}} \sum_{\substack{q \geq 17}} \sum_{\substack{9 < |a_q| \leq 2\sqrt{q} \\ (p+1-a_p) \mid (N+1-a_N) \\ (q+1-a_q) \mid (N+1-a_N) \\ 1-a_pa_q-p+pa_q \neq 0}}  \Pr\left[ \substack{ \#\Emod{p} = p +1 - a_p, \\ \#\Emod{q} = q+1-a_q} \biggr\vert p,q \text{ chosen} \right].
	\end{align*}
Corollary \ref{corProb} yields
	\begin{equation} \label{EqIneq5}
	\begin{aligned}
	&\sum_{\substack{q,a_q \\ q \geq 17 \\ 9 < |a_q| \leq 2\sqrt{q}}} \sum_{\substack{p,a_p \\ p \leq 13 \\ (p+1-a_p) \mid (N+1-a_N) \\ (q+1-a_q) \mid (N+1-a_N)}} \Pr\left[ \substack{ \#\Emod{p} = p +1 - a_p, \\ \#\Emod{q} = q+1-a_q} \biggr\vert p,q \text{ chosen} \right] \\
	&=\sum_{\substack{q,a_q \\ q \geq 17 \\ 9 < |a_q| \leq 2\sqrt{q}}} \sum_{\substack{p,a_p \\ p \leq 13 \\ (p+1-a_p) \mid (N+1-a_N) \\ (q+1-a_q) \mid (N+1-a_N) \\ 1-a_pa_q-p+pa_q = 0}} \Pr\left[ \substack{ \#\Emod{p} = p +1 - a_p, \\ \#\Emod{q} = q+1-a_q} \biggr\vert p,q \text{ chosen} \right] \\
	&+ \sum_{\substack{q,a_q \\ q \geq 17 \\ 9 < |a_q| \leq 2\sqrt{q}}} \sum_{\substack{p,a_p \\ p \leq 13 \\ (p+1-a_p) \mid (N+1-a_N) \\ (q+1-a_q) \mid (N+1-a_N) \\ 1-a_pa_q-p+pa_q \neq 0}} \Pr\left[ \substack{ \#\Emod{p} = p +1 - a_p, \\ \#\Emod{q} = q+1-a_q} \biggr\vert p,q \text{ chosen} \right] \\
	&\ll \sum_{\substack{p,a_p \\ p \leq 13}} \sum_{\substack{q \geq 17}} \sum_{\substack{9 < |a_q| \leq 2\sqrt{q} \\ (p+1-a_p) \mid (N+1-a_N) \\ (q+1-a_q) \mid (N+1-a_N) \\ 1-a_pa_q-p+pa_q = 0}}  \frac{1}{(pq)^{1/2-\epsilon}} \\
		&+ \sum_{\substack{p,a_p \\ p \leq 13}} \sum_{\substack{9 < |a_q| \leq 2\sqrt{M}}} \sum_{\substack{q \geq \left( \frac{a_q}{2} \right)^2 \\ (p+1-a_p) \mid (N+1-a_N) \\ (q+1-a_q) \mid (N+1-a_N) \\ 1-a_pa_q-p+pa_q \neq 0}}  \frac{1}{(pq)^{1/2-\epsilon}} \\
		&\ll \sum_{\substack{p,a_p \\ p \leq 13}} \sum_{\substack{q \geq 17}} \sum_{\substack{9 < |a_q| \leq 2\sqrt{q} \\ (p+1-a_p) \mid (N+1-a_N) \\ (q+1-a_q) \mid (N+1-a_N) \\ 1-a_pa_q-p+pa_q = 0}}  \frac{1}{q^{1/2-\epsilon}} \\
		&+ \sum_{\substack{p,a_p \\ p \leq 13}} \sum_{\substack{9 < |a_q| \leq 2\sqrt{M}}} \sum_{\substack{q \geq \left( \frac{a_q}{2} \right)^2 \\ (p+1-a_p) \mid (N+1-a_N) \\ (q+1-a_q) \mid (N+1-a_N) \\ 1-a_pa_q-p+pa_q \neq 0}}  \frac{1}{q^{1/2-\epsilon}}.
	\end{aligned}
	\end{equation}
Lemma \ref{LemmaFixPAP} bounds
	\begin{equation} \label{EqIneq6}
	\begin{aligned}
	&\sum_{\substack{p,a_p \\ p \leq 13}} \sum_{\substack{q \geq 17}} \sum_{\substack{9 < |a_q| \leq 2\sqrt{q} \\ (p+1-a_p) \mid (N+1-a_N) \\ (q+1-a_q) \mid (N+1-a_N) \\ 1-a_pa_q-p+pa_q = 0}}  \frac{1}{q^{1/2-\epsilon}} \\
	&\ll \sum_{\substack{p,a_p \\ p \leq 13}} \sum_{\substack{q \geq 17}}  \frac{1}{q^{1/2-\epsilon}} \\ 
	&\ll \sum_{\substack{p,a_p \\ p \leq 13}} \int_{17}^x \frac{1}{x^{1/2-\epsilon}} dx \\
	&\ll \sum_{\substack{p,a_p \\ p \leq 13}} M^{1/2+\epsilon} \\
	&\ll M^{1/2+\epsilon}
	\end{aligned}
	\end{equation}
and
	\begin{equation}\label{EqIneq7}
	\begin{aligned}
	&\sum_{\substack{p,a_p \\ p \leq 13}} \sum_{\substack{9 < |a_q| \leq 2\sqrt{M}}} \sum_{\substack{q \geq \left( \frac{a_q}{2} \right)^2 \\ (p+1-a_p) \mid (N+1-a_N) \\ (q+1-a_q) \mid (N+1-a_N) \\ 1-a_pa_q-p+pa_q \neq 0}}  \frac{1}{q^{1/2-\epsilon}} \\
	&\ll \sum_{\substack{p,a_p \\ p \leq 13}} \sum_{\substack{9 < |a_q| \leq 2\sqrt{M}}} \sum_{\substack{q \geq \left( \frac{a_q}{2} \right)^2 \\ (p+1-a_p) \mid (N+1-a_N) \\ (q+1-a_q) \mid (N+1-a_N) \\ 1-a_pa_q-p+pa_q \neq 0}}  \frac{1}{a_q^{1-\epsilon}} \\
	&\ll \sum_{\substack{p,a_p \\ p \leq 13}} \sum_{\substack{9 < |a_q| \leq 2\sqrt{M}}} \frac{1}{a_q^{1-\epsilon}} \\
	&\ll \sum_{\substack{p,a_p \\ p \leq 13}} \int_9^{2\sqrt{M}} \frac{1}{x^{1-\epsilon}} dx \\
	&\ll \sum_{\substack{p,a_p \\ p \leq 13}} M^{\epsilon/2} \\
	&\ll M^{\epsilon/2}.
	\end{aligned}
	\end{equation}
Combining (\ref{EqIneq5}), (\ref{EqIneq6}), and (\ref{EqIneq7}) we have that
	\begin{equation}\label{EqIneq8}
	\begin{aligned}
		&\sum_{\substack{q,a_q \\ q \geq 17 \\ 9 < |a_q| \leq 2\sqrt{q}}} \sum_{\substack{p,a_p \\ p \leq 13 \\ (p+1-a_p) \mid (N+1-a_N) \\ (q+1-a_q) \mid (N+1-a_N)}} \Pr\left[ \substack{ \#\Emod{p} = p +1 - a_p, \\ \#\Emod{q} = q+1-a_q} \biggr\vert p,q \text{ chosen} \right]  \\
		&\ll M^{1/2+\epsilon}.
	\end{aligned}
	\end{equation}
Corollary \ref{corProb} bounds the second sum of (\ref{EqIneq10}) as
	\begin{align*}
	&\sum_{\substack{q,a_q \\ q \geq 17 \\ 9 < |a_q| \leq 2\sqrt{q}}} \sum_{\substack{p,a_p \\ \frac{\sqrt{q}}{16} \leq p < q \\ (p+1-a_p) \mid (N+1-a_N) \\ (q+1-a_q) \mid (N+1-a_N)}} \Pr\left[ \substack{ \#\Emod{p} = p +1 - a_p, \\ \#\Emod{q} = q+1-a_q} \biggr\vert p,q \text{ chosen} \right] \\
	&\ll\sum_{\substack{q,a_q \\ q \geq 17 \\ 9 < |a_q| \leq 2\sqrt{q}}} \sum_{\substack{p,a_p \\ \frac{\sqrt{q}}{16} \leq p < q \\ (p+1-a_p) \mid (N+1-a_N) \\ (q+1-a_q) \mid (N+1-a_N)}} \frac{(4q - a_q^2)^{1/2+\epsilon}}{q^{1+\epsilon} p^{1/2-\epsilon}} \\
	&\ll \sum_{\substack{q,a_q \\ q \geq 17 \\ 9 < |a_q| \leq 2\sqrt{q}}} \sum_{\substack{p,a_p \\ \frac{\sqrt{q}}{16} \leq p < q \\ (p+1-a_p) \mid (N+1-a_N) \\ (q+1-a_q) \mid (N+1-a_N)}} \frac{(4q - a_q^2)^{1/2+\epsilon}}{q^{1+\epsilon} (q^{1/2})^{1/2-\epsilon}} \\
	&= \sum_{\substack{q,a_q \\ q \geq 17 \\ 9 < |a_q| \leq 2\sqrt{q}}} \frac{1}{q^{5/4+\epsilon/2}} \sum_{\substack{p,a_p \\ \frac{\sqrt{q}}{16} \leq p < q \\ (p+1-a_p) \mid (N+1-a_N) \\ (q+1-a_q) \mid (N+1-a_N)}} (4q - a_q^2)^{1/2+\epsilon}.
	\end{align*}
Lemma \ref{LemmaBoundPair1} and Lemma \ref{LemmaBoundPair2} show that each choice of $q$ and $a_q$ in the above sum yield $O(a_qq^\epsilon)$ possible choices of $p$ and $a_p$. Thus,
	\begin{align*}
	&\sum_{\substack{q,a_q \\ q \geq 17 \\ 9 < |a_q| \leq 2\sqrt{q}}} \sum_{\substack{p,a_p \\ \frac{\sqrt{q}}{16} \leq p < q \\ (p+1-a_p) \mid (N+1-a_N) \\ (q+1-a_q) \mid (N+1-a_N)}} \Pr\left[ \substack{ \#\Emod{p} = p +1 - a_p, \\ \#\Emod{q} = q+1-a_q} \biggr\vert p,q \text{ chosen} \right] \\
	&\ll \sum_{\substack{q,a_q \\ q \geq 17 \\ 9 < |a_q| \leq 2\sqrt{q}}} \frac{1}{q^{5/4+\epsilon/2}} a_qq^{\epsilon} (4q - a_q^2)^{1/2+\epsilon} \\
	&\ll \sum_{\substack{q,a_q \\ q \geq 17 \\ 9 < |a_q| \leq 2\sqrt{q}}} \frac{1}{q^{5/4+\epsilon/2}} a_qq^{\epsilon} (4q - a_q^2)^{1/2} q^{\epsilon} \\
	&=\sum_{\substack{q,a_q \\ q \geq 17 \\ 9 < |a_q| \leq 2\sqrt{q}}} \frac{1}{q^{5/4-3\epsilon/2}} a_q (4q - a_q^2)^{1/2} \\
	&= \sum_{q \geq 17} \frac{1}{q^{5/4-3\epsilon/2}} \sum_{9 < |a_q| \leq 2\sqrt{q}} a_q (4q - a_q^2)^{1/2} \\
	&= \sum_{q \geq 17} \frac{1}{q^{5/4-3\epsilon/2}} \sum_{9 < |a_q| \leq 2\sqrt{q}} 4q \frac{a_q}{2\sqrt{q}} \left(1 - \left( \frac{a_q}{2\sqrt{q}} \right)^2 \right)^{1/2} \\
	&\ll \sum_{q \geq 17} \frac{1}{q^{1/4-3\epsilon/2}} \sum_{9 < |a_q| \leq 2\sqrt{q}} \frac{a_q}{2\sqrt{q}} \left(1 - \left( \frac{a_q}{2\sqrt{q}} \right)^2 \right)^{1/2}.
	\end{align*}
Note that
	\begin{align*}
		\sum_{9 < |a_q| \leq 2\sqrt{q}} \frac{a_q}{2\sqrt{q}} \left(1 - \left( \frac{a_q}{2\sqrt{q}} \right)^2 \right)^{1/2} \ll \int_{0}^1 x \sqrt{1-x^2} dx = O(1).
	\end{align*}
Therefore,
	\begin{equation} \label{EqIneq9}
	\begin{aligned}
		&\sum_{\substack{q,a_q \\ q \geq 17 \\ 9 < |a_q| \leq 2\sqrt{q}}} \sum_{\substack{p,a_p \\ \frac{\sqrt{q}}{16} \leq p < q \\ (p+1-a_p) \mid (N+1-a_N) \\ (q+1-a_q) \mid (N+1-a_N)}} \Pr\left[ \substack{ \#\Emod{p} = p +1 - a_p, \\ \#\Emod{q} = q+1-a_q} \biggr\vert p,q \text{ chosen} \right] \\
		&\ll \sum_{q \geq 17} \frac{1}{q^{1/4-\epsilon/2}} \\
		&\ll \int_{17}^M \frac{1}{x^{1/4-3\epsilon/2}} \\
		&\ll M^{3/4+3\epsilon/2}.
	\end{aligned}
	\end{equation} 
Combining (\ref{EqIneq10}), (\ref{EqIneq8}), and (\ref{EqIneq9}) and replacing $\epsilon$ with $2\epsilon/3$ yields
	\begin{equation} \label{EqIneq11}
	\begin{aligned}
	&\sum_{\substack{q,a_q \\ q \geq 17 \\ 9 < |a_q| \leq 2\sqrt{q}}} \sum_{\substack{p,a_p \\ p < q \\ (p+1-a_p) \mid (N+1-a_N) \\ (q+1-a_q) \mid (N+1-a_N)}} \Pr\left[ \substack{ \#\Emod{p} = p +1 - a_p, \\ \#\Emod{q} = q+1-a_q} \biggr\vert p,q \text{ chosen} \right]\\
	&\ll M^{3/4+\epsilon}.
	\end{aligned}
	\end{equation}
Furthermore, (\ref{EqIneq12}), (\ref{EqIneq13}), (\ref{EqIneq4}), and (\ref{EqIneq11}) altogether bound
	\begin{align*}
		\sum_{\substack{q,a_q,p,a_p \\ (p+1-a_p),(q+1-a_q) \mid (N+1-a_N)  \\ a_q \neq 1 \\ p < q} } \Pr\left[ \substack{ \#\Emod{p} = p +1 - a_p, \\ \#\Emod{q} = q+1-a_q} \biggr\vert p,q \text{ chosen} \right] \ll M^{3/4+\epsilon}.
	\end{align*}
Going back to (\ref{EqIneq14}) we have that
	\begin{align*}
		 &\Pr \left[\substack{(p+1-a_p),(q+1-a_q) \mid (N+1-a_N), \\
					a_p \text{ or } a_q \neq 1} \right] \\
		&\ll \left( \frac{\log M}{M} \right)^2 M^{3/4+\epsilon} \\
		&= \frac{ M^\epsilon (\log M)^2}{M^{5/4}} \\
		&\ll \frac{1}{M^{5/4-\epsilon'}}
	\end{align*}
	for all $\epsilon' > 0$ as desired.
\end{proof}
\end{lem}
Again, Corollary \ref{corEnd} proves the conjecture stated in \cite{REU2016}.
\begin{cor} \label{corEnd}
Let $5 \leq p,q \leq M$ be randomly chosen distinct primes and let $N = pq$. Let $\Emod{N}$ be a randomly chosen elliptic curve with good reduction at $p$ and $q$ such that $(p+1-a_p), (q+1-a_q) \mid (N+1-a_N)$. Then\begin{align*}
	\lim_{M \rightarrow \infty} \Pr[a_p \text{ or } a_q \text{ is not 1}] = 0
\end{align*} 
and
\begin{align*}
	\lim_{M \rightarrow \infty} \Pr[\#\Emod{N} = N+1-a_N] = 1.
\end{align*}
\begin{proof}
When $E$ is a random elliptic curve with good reduction at $p$ and $q$, not necessarily with $(p+1-a_p),(q+1-a_q) \mid (N+1-a_N)$, Lemma \ref{LemmaLBound} and Lemma \ref{LemmaUBound} show that
	\begin{align*}
		\frac{ \Pr[a_p \text{ or } a_q \text{ is not 1 and } (p+1-a_p), (q+1-a_q) \mid N+1-a_N]}{\Pr[a_p = a_q = 1]} &\ll \frac{\frac{1}{M^{{5/4}-\epsilon}}}{\frac{1}{M^{1+\epsilon}}} \\
		&\ll \frac{1}{M^{1/4-2\epsilon}}.
	\end{align*}
Thus, if $E$ satisfies $(p+1-a_p), (q+1-a_q) \mid (N+1-a_N)$ we have 
	\begin{align*}
	\lim_{M \rightarrow \infty} \Pr[a_p \text{ or } a_q \text{ is not 1}] = 0.
	\end{align*} 
Since $\#\Emod{N} \neq N+1-a_N$ we have that $a_p$ or $a_q$ is not $1$,
	\begin{align*}
		\lim_{M \rightarrow \infty} \Pr[\#\Emod{N} \neq N+1-a_N] = 0,
	\end{align*}
we have that $\lim_{M \rightarrow \infty} \Pr[\#\Emod{N} = N+1-a_N] = 1$.
	
} 
\else{
Fix $M \geq 7$. Let $p$ and $q$ be primes with $5 \leq p,q \leq M$, let $a_p$ and $a_q$ be integers such that $|a_p| \leq 2\sqrt{p}$ and $|a_q| \leq 2\sqrt{q}$, and let $a_N = a_pa_q$. Let $T$ be the set
	\begin{align*}
		T = \left \{ (q,a_q,p,a_p) \in \mathbb{Z}^4 \bigg\vert \begin{gathered} p,q \text{ prime, } 5 \leq p,q \leq M, \; |a_p| \leq 2\sqrt{p}, \; |a_q| \leq 2\sqrt{q}, \\ a_p \text{ or } a_q \neq 1, \; (p+1-a_p),(q+1-a_q) \mid (N+1-a_N)\end{gathered} \right \}.
	\end{align*}
Furthermore, let $S = \left \{ (q,a_q,p,a_p) \in T \mid p < q \right \}$.

By Lemma \ref{Lemmaapaq1}, $a_q\neq1$ for every $(q,a_q,p,a_p) \in S$.  The number of pairs of $p$ and $q$ is on the order of  $\left( \frac{M}{\log M} \right)^2$. Note that
	\begin{equation} \label{EqIneq14}
	\begin{aligned}
		 &\Pr \left[\substack{(p+1-a_p),(q+1-a_q) \mid (N+1-a_N), \\
					a_p \text{ or } a_q \neq 1} \right] \\
		&= \sum_{(q,a_q,p,a_p) \in T } \Pr\left[ \substack{ p,q \text{ chosen}, \\ \#\Emod{p} = p +1 - a_p, \\ \#\Emod{q} = q+1-a_q} \right] \\
		&= \sum_{(q,a_q,p,a_p) \in T } \Pr\left[ \substack{ \#\Emod{p} = p +1 - a_p, \\ \#\Emod{q} = q+1-a_q} \biggr\vert p,q \text{ chosen} \right]  \\
				&\phantom{=} \cdot \left( \frac{1}{\#\{(p,q) \mid p,q \text{ distinct primes}, 5 \leq p,q \leq M\}} \right) \\
		&\approx \left( \frac{\log M}{M} \right)^2 \sum_{(q,a_q,p,a_p) \in T } \Pr\left[ \substack{ \#\Emod{p} = p +1 - a_p, \\ \#\Emod{q} = q+1-a_q} \biggr\vert p,q \text{ chosen} \right].
	\end{aligned}
	\end{equation}
Furthermore,
	\begin{align*}
		 &\sum_{(q,a_q,p,a_p) \in T } \Pr\left[ \substack{ \#\Emod{p} = p +1 - a_p, \\ \#\Emod{q} = q+1-a_q} \biggr\vert p,q \text{ chosen} \right] \\
		&=2 \cdot \left(\sum_{\substack{(q,a_q,p,a_p) \in T \\ p < q}} \Pr\left[ \substack{ \#\Emod{p} = p +1 - a_p, \\ \#\Emod{q} = q+1-a_q} \biggr\vert p,q \text{ chosen} \right] \right) \\
		&\approx \sum_{\substack{(q,a_q,p,a_p) \in T \\ p < q}} \Pr\left[ \substack{ \#\Emod{p} = p +1 - a_p, \\ \#\Emod{q} = q+1-a_q} \biggr\vert p,q \text{ chosen} \right] \\ 
		&= \sum_{(q,a_q,p,a_p) \in S} \Pr\left[ \substack{ \#\Emod{p} = p +1 - a_p, \\ \#\Emod{q} = q+1-a_q} \biggr\vert p,q \text{ chosen} \right].
\end{align*}
By Lemma \ref{Lemmaapaq1},
\begin{align*}
	&\sum_{(q,a_q,p,a_p) \in S} \Pr\left[ \substack{ \#\Emod{p} = p +1 - a_p, \\ \#\Emod{q} = q+1-a_q} \biggr\vert p,q \text{ chosen} \right] =1.\\
\end{align*}
Moreover,
\begin{equation}\label{EqIneq12}
\begin{aligned}
		&\sum_{(q,a_q,p,a_p) \in S} \Pr\left[ \substack{ \#\Emod{p} = p +1 - a_p, \\ \#\Emod{q} = q+1-a_q} \biggr\vert p,q \text{ chosen} \right] \\
		&= \sum_{\substack{q,a_q \\ q \geq 17 \\ 9 < |a_q| \leq 2\sqrt{q}}} \sum_{(q,a_q,p,a_p) \in S} \Pr\left[ \substack{ \#\Emod{p} = p +1 - a_p, \\ \#\Emod{q} = q+1-a_q} \biggr\vert p,q \text{ chosen} \right] \\
		&+ \sum_{q \geq 17} \sum_{\substack{ |a_q| \leq 9  \\ a_q \neq 1} } \sum_{(q,a_q,p,a_p) \in S} \Pr\left[ \substack{ \#\Emod{p} = p +1 - a_p, \\ \#\Emod{q} = q+1-a_q} \biggr\vert p,q \text{ chosen} \right] \\
		&+ \sum_{\substack{(q,a_q,p,a_p) \in S \\ q < 17}}\Pr\left[ \substack{ \#\Emod{p} = p +1 - a_p, \\ \#\Emod{q} = q+1-a_q} \biggr\vert p,q \text{ chosen} \right].
	\end{aligned}
\end{equation}
Next, we bound the three parts to the above sum. Note that if $p < q < 17$, then there are only $O(1)$ possible combinations of $q,a_q,p$ and $a_p$. Thus,
	\begin{equation}\label{EqIneq13}
	\begin{aligned}
		\sum_{\substack{(q,a_q,p,a_p) \in S \\ q < 17}}\Pr\left[ \substack{ \#\Emod{p} = p +1 - a_p, \\ \#\Emod{q} = q+1-a_q} \biggr\vert p,q \text{ chosen} \right]. 
	\end{aligned}
	\end{equation}
By Lemma \ref{LemmaAnomProd}, 
	\begin{equation}\label{EqIneq1}
	\begin{aligned} 
		&\sum_{q \geq 17} \sum_{\substack{ |a_q| \leq 9  \\ a_q \neq 1} } \sum_{(q,a_q,p,a_p) \in S} \Pr\left[ \substack{ \#\Emod{p} = p +1 - a_p, \\ \#\Emod{q} = q+1-a_q} \biggr\vert p,q \text{ chosen} \right] \\
		&= \sum_{q \geq 17} \sum_{\substack{ |a_q| \leq 9  \\ a_q \neq 1} } \sum_{\substack{(q,a_q,p,a_p) \in S \\ p \leq 13}} \Pr\left[ \substack{ \#\Emod{p} = p +1 - a_p, \\ \#\Emod{q} = q+1-a_q} \biggr\vert p,q \text{ chosen} \right] \\
		&+ \sum_{q \geq 17} \sum_{\substack{ |a_q| \leq 9  \\ a_q \neq 1} } \sum_{\substack{(q,a_q,p,a_p) \in S \\ 13 < p}} \Pr\left[ \substack{ \#\Emod{p} = p +1 - a_p, \\ \#\Emod{q} = q+1-a_q} \biggr\vert p,q \text{ chosen} \right].
	\end{aligned}
	\end{equation}
By Corollary \ref{corProb} the first summand of (\ref{EqIneq1}) satisfies
	\begin{align*}
		&\sum_{q \geq 17} \sum_{\substack{ |a_q| \leq 9  \\ a_q \neq 1} } \sum_{\substack{(q,a_q,p,a_p) \in S \\ p \leq 13}} \Pr\left[ \substack{ \#\Emod{p} = p +1 - a_p, \\ \#\Emod{q} = q+1-a_q} \biggr\vert p,q \text{ chosen} \right] \\
		&\ll \sum_{q \geq 17} \sum_{\substack{ |a_q| \leq 9  \\ a_q \neq 1} } \sum_{\substack{(q,a_q,p,a_p) \in S \\ p \leq 13}} \frac{1}{q^{1/2-\epsilon}}.
	\end{align*}
There are $O(1)$ possible combinations of $a_q,p$ and $a_p$ in the above sum, so
	\begin{equation} \label{EqIneq2}
	\begin{aligned}
		&\sum_{q \geq 17} \sum_{\substack{ |a_q| \leq 9  \\ a_q \neq 1} } \sum_{\substack{(q,a_q,p,a_p) \in S \\ p \leq 13}} \Pr\left[ \substack{ \#\Emod{p} = p +1 - a_p, \\ \#\Emod{q} = q+1-a_q} \biggr\vert p,q \text{ chosen} \right] \\
		&\ll \sum_{q \geq 17} \frac{1}{q^{1/2-\epsilon}} 
		\ll \int_{17}^{M} \frac{1}{x^{1/2-\epsilon}} dx
		\ll M^{1/2+\epsilon}.
	\end{aligned}
	\end{equation}
On the other hand, $|a_p| \leq 2\sqrt{p} < 2\sqrt{q}$, so Lemma \ref{LemmaNumAp} and Lemma \ref{LemmaBounds} show that the second summand of (\ref{EqIneq1}) satisfies
	\begin{align*}
		 &\sum_{q \geq 17} \sum_{\substack{ |a_q| \leq 9  \\ a_q \neq 1} } \sum_{\substack{(q,a_q,p,a_p) \in S \\ 13 < p}} \Pr\left[ \substack{ \#\Emod{p} = p +1 - a_p, \\ \#\Emod{q} = q+1-a_q} \biggr\vert p,q \text{ chosen} \right] \\
		&\ll  \sum_{q \geq 17} \sum_{\substack{ |a_q| \leq 9  \\ a_q \neq 1} } \sum_{|a_p| \leq 2\sqrt{q}} \sum_{\substack{(q,a_q,p,a_p) \in S \\ 13 < p}} \Pr\left[ \substack{ \#\Emod{p} = p +1 - a_p, \\ \#\Emod{q} = q+1-a_q} \biggr\vert p,q \text{ chosen} \right] \\
		&\ll \sum_{q \geq 17} \sum_{\substack{ |a_q| \leq 9  \\ a_q \neq 1} } \sum_{|a_p| \leq 2\sqrt{q}} \sum_{\substack{(q,a_q,p,a_p) \in S \\ 13 < p}} \frac{1}{(pq)^{1/2-\epsilon}}.
	\end{align*}
By Lemma \ref{LemmaAnomProd}, $p \geq \frac{\sqrt{q}}{16}$ in the above sum. Thus,
	\begin{align*}
		 &\sum_{q \geq 17} \sum_{\substack{ |a_q| \leq 9  \\ a_q \neq 1} } \sum_{\substack{(q,a_q,p,a_p) \in S \\ 13 < p}} \Pr\left[ \substack{ \#\Emod{p} = p +1 - a_p, \\ \#\Emod{q} = q+1-a_q} \biggr\vert p,q \text{ chosen} \right] \\
		&\ll \sum_{q \geq 17} \sum_{\substack{ |a_q| \leq 9  \\ a_q \neq 1} } \sum_{|a_p| \leq 2\sqrt{q}} \sum_{\substack{(q,a_q,p,a_p) \in S \\ \frac{\sqrt{q}}{16} \leq p}} \frac{1}{(pq)^{1/2-\epsilon}} \\
		&= \sum_{q \geq 17} \sum_{\substack{ |a_q| \leq 9  \\ a_q \neq 1} } \sum_{|a_p| \leq 2\sqrt{q}} \sum_{\substack{(q,a_q,p,a_p) \in S \\ \frac{\sqrt{q}}{16} \leq p}} \frac{1}{q^{3/4-3\epsilon/2}}.
	\end{align*}
Lemma \ref{LemmaFixQAQ} implies that each combination of $q,a_q$ and $a_p$ yields only $O(1)$ possible values of $p$, so
	\begin{align*}
	&\sum_{q \geq 17} \sum_{\substack{ |a_q| \leq 9  \\ a_q \neq 1} } \sum_{\substack{(q,a_q,p,a_p) \in S \\ 13 < p}} \Pr\left[ \substack{ \#\Emod{p} = p +1 - a_p, \\ \#\Emod{q} = q+1-a_q} \biggr\vert p,q \text{ chosen} \right] \\
	&\ll \sum_{q \geq 17} \sum_{\substack{ |a_q| \leq 9  \\ a_q \neq 1} } \sum_{|a_p| \leq 2\sqrt{q}} \frac{1}{q^{3/4-3\epsilon/2}}  \\
	&\ll \int_{17}^{M} \frac{1}{x^{1/4-3\epsilon/2}} dx 
	\ll M^{3/4+3\epsilon/2}.
	\end{align*}
Replacing $\epsilon$ with $2\epsilon/3$, we have 
	\begin{equation} \label{EqIneq3}
	\begin{aligned}
		&\sum_{q \geq 17} \sum_{\substack{ |a_q| \leq 9  \\ a_q \neq 1} } \sum_{\substack{(q,a_q,p,a_p) \in S \\ 13 < p}} \Pr\left[ \substack{ \#\Emod{p} = p +1 - a_p, \\ \#\Emod{q} = q+1-a_q} \biggr\vert p,q \text{ chosen} \right] \\\
		&\ll M^{3/4+\epsilon}
	\end{aligned}
	\end{equation}
for all $\epsilon > 0$. Combining (\ref{EqIneq1}), (\ref{EqIneq2}) and (\ref{EqIneq3}) shows that
	\begin{equation} \label{EqIneq4}
	\begin{aligned}
		&\sum_{q \geq 17} \sum_{\substack{ |a_q| \leq 9  \\ a_q \neq 1} } \sum_{(q,a_q,p,a_p) \in S} \Pr\left[ \substack{ \#\Emod{p} = p +1 - a_p, \\ \#\Emod{q} = q+1-a_q} \biggr\vert p,q \text{ chosen} \right] \\
		&\ll M^{3/4+\epsilon}.
	\end{aligned}
	\end{equation}
Using Lemma \ref{LemmaAnomProd}, express the first summand of (\ref{EqIneq12}) as 
	\begin{equation} \label{EqIneq10}
	\begin{aligned}
		& \sum_{\substack{q,a_q \\ q \geq 17 \\ 9 < |a_q| \leq 2\sqrt{q}}} \sum_{(q,a_q,p,a_p) \in S} \Pr\left[ \substack{ \#\Emod{p} = p +1 - a_p, \\ \#\Emod{q} = q+1-a_q} \biggr\vert p,q \text{ chosen} \right]\\ 
		&=  \sum_{\substack{q,a_q \\ q \geq 17 \\ 9 < |a_q| \leq 2\sqrt{q}}} \sum_{\substack{(q,a_q,p,a_p) \in S \\ p \leq 13}} \Pr\left[ \substack{ \#\Emod{p} = p +1 - a_p, \\ \#\Emod{q} = q+1-a_q} \biggr\vert p,q \text{ chosen} \right] \\
		&+ \sum_{\substack{q,a_q \\ q \geq 17 \\ 9 < |a_q| \leq 2\sqrt{q}}} \sum_{\substack{(q,a_q,p,a_p) \in S \\ \frac{\sqrt{q}}{16} \leq p}} \Pr\left[ \substack{ \#\Emod{p} = p +1 - a_p, \\ \#\Emod{q} = q+1-a_q} \biggr\vert p,q \text{ chosen} \right].
	\end{aligned}
	\end{equation}	
By rearranging the first of the two summands of (\ref{EqIneq10}) we obtain
	\begin{align*}
		&\sum_{\substack{q,a_q \\ q \geq 17 \\ 9 < |a_q| \leq 2\sqrt{q}}} \sum_{\substack{(q,a_q,p,a_p) \in S \\ p \leq 13}} \Pr\left[ \substack{ \#\Emod{p} = p +1 - a_p, \\ \#\Emod{q} = q+1-a_q} \biggr\vert p,q \text{ chosen} \right] \\
                 &= \sum_{\substack{p,a_p \\ p \leq 13}} \sum_{\substack{q \geq 17}} \sum_{\substack{(q,a_q,p,a_p) \in S \\ 9 < |a_q| \\ 1-a_pa_q - p + pa_q = 0} }  \Pr\left[ \substack{ \#\Emod{p} = p +1 - a_p, \\ \#\Emod{q} = q+1-a_q} \biggr\vert p,q \text{ chosen} \right] \\
		&+ \sum_{\substack{p,a_p \\ p \leq 13}} \sum_{\substack{q \geq 17}} \sum_{\substack{(q,a_q,p,a_p) \in S \\ 9 < |a_q| \\ 1-a_pa_q - p + pa_q \neq 0} }  \Pr\left[ \substack{ \#\Emod{p} = p +1 - a_p, \\ \#\Emod{q} = q+1-a_q} \biggr\vert p,q \text{ chosen} \right].
	\end{align*}
Corollary \ref{corProb} yields
	\begin{equation} \label{EqIneq5}
	\begin{aligned}
	&\sum_{\substack{q,a_q \\ q \geq 17 \\ 9 < |a_q| \leq 2\sqrt{q}}} \sum_{\substack{(q,a_q,p,a_p) \in S \\ p \leq 13}} \Pr\left[ \substack{ \#\Emod{p} = p +1 - a_p, \\ \#\Emod{q} = q+1-a_q} \biggr\vert p,q \text{ chosen} \right] \\
	&\ll \sum_{\substack{p,a_p \\ p \leq 13}} \sum_{\substack{q \geq 17}} \sum_{\substack{(q,a_q,p,a_p) \in S \\ 9 < |a_q| \\ 1-a_pa_q - p + pa_q = 0} }  \frac{1}{q^{1/2-\epsilon}} 
	+ \sum_{\substack{p,a_p \\ p \leq 13}} \sum_{\substack{9 < |a_q| \leq 2\sqrt{M}}} \sum_{\substack{(q,a_q,p,a_p) \in S \\ q \geq \left( \frac{a_q}{2} \right)^2 \\ 1-a_pa_q - p + pa_q \neq 0} }  \frac{1}{q^{1/2-\epsilon}}.
	\end{aligned}
	\end{equation}
Lemma \ref{LemmaFixPAP} bounds
	\begin{equation} \label{EqIneq6}
	\begin{aligned}
	&\sum_{\substack{p,a_p \\ p \leq 13}} \sum_{\substack{q \geq 17}} \sum_{\substack{(q,a_q,p,a_p) \in S \\ 9 < |a_q| \\ 1-a_pa_q - p + pa_q = 0} }  \frac{1}{q^{1/2-\epsilon}}
	\ll \sum_{\substack{p,a_p \\ p \leq 13}} \sum_{\substack{q \geq 17}}  \frac{1}{q^{1/2-\epsilon}} 
	\ll \sum_{\substack{p,a_p \\ p \leq 13}} M^{1/2+\epsilon} 
	\ll M^{1/2+\epsilon}
	\end{aligned}
	\end{equation}
and
	\begin{equation}\label{EqIneq7}
	\begin{aligned}
	&\sum_{\substack{p,a_p \\ p \leq 13}} \sum_{\substack{9 < |a_q| \leq 2\sqrt{M}}} \sum_{\substack{(q,a_q,p,a_p) \in S \\ q \geq \left( \frac{a_q}{2} \right)^2 \\ 1-a_pa_q - p + pa_q \neq 0} }  \frac{1}{q^{1/2-\epsilon}} \\
	&\ll \sum_{\substack{p,a_p \\ p \leq 13}} \sum_{\substack{9 < |a_q| \leq 2\sqrt{M}}} \frac{1}{a_q^{1-2\epsilon}} \\
	&\ll \sum_{\substack{p,a_p \\ p \leq 13}} \int_9^{2\sqrt{M}} \frac{1}{x^{1-2\epsilon}} dx 
	\ll \sum_{\substack{p,a_p \\ p \leq 13}} M^{\epsilon} 
	\ll M^{\epsilon}.
	\end{aligned}
	\end{equation}
Combining (\ref{EqIneq5}), (\ref{EqIneq6}), and (\ref{EqIneq7}) we have that
	\begin{equation}\label{EqIneq8}
		\sum_{\substack{q,a_q \\ q \geq 17 \\ 9 < |a_q| \leq 2\sqrt{q}}} \sum_{\substack{(q,a_q,p,a_p) \in S \\ p \leq 13}} \Pr\left[ \substack{ \#\Emod{p} = p +1 - a_p, \\ \#\Emod{q} = q+1-a_q} \biggr\vert p,q \text{ chosen} \right] 
		\ll M^{1/2+\epsilon}
	\end{equation}
On the other hand, Corollary \ref{corProb} bounds the second sum of (\ref{EqIneq10}) as
	\begin{align*}
	& \sum_{\substack{q,a_q \\ q \geq 17 \\ 9 < |a_q| \leq 2\sqrt{q}}} \sum_{\substack{(q,a_q,p,a_p) \in S \\ \frac{\sqrt{q}}{16} \leq p}} \Pr\left[ \substack{ \#\Emod{p} = p +1 - a_p, \\ \#\Emod{q} = q+1-a_q} \biggr\vert p,q \text{ chosen} \right] \\
	&\ll \sum_{\substack{q,a_q \\ q \geq 17 \\ 9 < |a_q| \leq 2\sqrt{q}}} \sum_{\substack{(q,a_q,p,a_p) \in S \\ \frac{\sqrt{q}}{16} \leq p}} \frac{(4q - a_q^2)^{1/2+\epsilon}}{q^{1+\epsilon} p^{1/2-\epsilon}} \\
	&=  \sum_{\substack{q,a_q \\ q \geq 17 \\ 9 < |a_q| \leq 2\sqrt{q}}} \frac{1}{q^{5/4+\epsilon/2}} \sum_{\substack{(q,a_q,p,a_p) \in S \\ \frac{\sqrt{q}}{16} \leq p}} (4q - a_q^2)^{1/2+\epsilon}.
	\end{align*}
	Lemmas \ref{LemmaBoundPair1} and \ref{LemmaBoundPair2} show that each choice of $q$ and $a_q$ in the above sum yield $O(a_qq^\epsilon)$ possible choices of $p$ and $a_p$. Thus,
	\begin{align*}
	&\sum_{\substack{q,a_q \\ q \geq 17 \\ 9 < |a_q| \leq 2\sqrt{q}}} \sum_{\substack{(q,a_q,p,a_p) \in S \\ \frac{\sqrt{q}}{16} \leq p}} \Pr\left[ \substack{ \#\Emod{p} = p +1 - a_p, \\ \#\Emod{q} = q+1-a_q} \biggr\vert p,q \text{ chosen} \right] \\
	&\ll \sum_{\substack{q,a_q \\ q \geq 17 \\ 9 < |a_q| \leq 2\sqrt{q}}} \frac{1}{q^{5/4+\epsilon/2}} a_qq^{\epsilon} (4q - a_q^2)^{1/2+\epsilon} \\
	&\ll \sum_{\substack{q,a_q \\ q \geq 17 \\ 9 < |a_q| \leq 2\sqrt{q}}} \frac{1}{q^{5/4+\epsilon/2}} a_qq^{\epsilon} (4q - a_q^2)^{1/2} q^{\epsilon} \\
	&= \sum_{q \geq 17} \frac{1}{q^{5/4-3\epsilon/2}} \sum_{9 < |a_q| \leq 2\sqrt{q}} 4q \frac{a_q}{2\sqrt{q}} \left(1 - \left( \frac{a_q}{2\sqrt{q}} \right)^2 \right)^{1/2} \\
	&\ll \sum_{q \geq 17} \frac{1}{q^{1/4-3\epsilon/2}} \sum_{9 < |a_q| \leq 2\sqrt{q}} \frac{a_q}{2\sqrt{q}} \left(1 - \left( \frac{a_q}{2\sqrt{q}} \right)^2 \right)^{1/2}.
	\end{align*}	
Note that
	\begin{align*}
		\sum_{9 < |a_q| \leq 2\sqrt{q}} \frac{a_q}{2\sqrt{q}} \left(1 - \left( \frac{a_q}{2\sqrt{q}} \right)^2 \right)^{1/2} \ll \int_{0}^1 x \sqrt{1-x^2} dx = O(1).
	\end{align*}
Therefore,
	\begin{equation} \label{EqIneq9}
	\begin{aligned}
		&\sum_{\substack{q,a_q \\ q \geq 17 \\ 9 < |a_q| \leq 2\sqrt{q}}} \sum_{\substack{(q,a_q,p,a_p) \in S \\ \frac{\sqrt{q}}{16} \leq p}} \Pr\left[ \substack{ \#\Emod{p} = p +1 - a_p, \\ \#\Emod{q} = q+1-a_q} \biggr\vert p,q \text{ chosen} \right] \\
		&\ll \sum_{q \geq 17} \frac{1}{q^{1/4-\epsilon/2}} 
		\ll \int_{17}^M \frac{1}{x^{1/4-3\epsilon/2}} 
		\ll M^{3/4+3\epsilon/2}.
	\end{aligned}
	\end{equation} 
Combining (\ref{EqIneq10}), (\ref{EqIneq8}), and (\ref{EqIneq9}) and replacing $\epsilon$ with $2\epsilon/3$ yields
	\begin{equation} \label{EqIneq11}
	\begin{aligned}
	& \sum_{\substack{q,a_q \\ q \geq 17 \\ 9 < |a_q| \leq 2\sqrt{q}}} \sum_{(q,a_q,p,a_p) \in S} \Pr\left[ \substack{ \#\Emod{p} = p +1 - a_p, \\ \#\Emod{q} = q+1-a_q} \biggr\vert p,q \text{ chosen} \right]
	\ll M^{3/4+\epsilon}.
	\end{aligned}
	\end{equation}
Furthermore, (\ref{EqIneq12}), (\ref{EqIneq13}), (\ref{EqIneq4}), and (\ref{EqIneq11}) altogether bound
	\begin{align*}
		\sum_{(q,a_q,p,a_p) \in S} \Pr\left[ \substack{ \#\Emod{p} = p +1 - a_p, \\ \#\Emod{q} = q+1-a_q} \biggr\vert p,q \text{ chosen} \right] \ll M^{3/4+\epsilon}.
	\end{align*}
From (\ref{EqIneq14}) we have 
	\begin{align*}
		 &\Pr \left[\substack{(p+1-a_p),(q+1-a_q) \mid (N+1-a_N), \\
					a_p \text{ or } a_q \neq 1} \right] \\
		&\ll \left( \frac{\log M}{M} \right)^2 M^{3/4+\epsilon} 
		= \frac{ M^\epsilon (\log M)^2}{M^{5/4}} 
		\ll \frac{1}{M^{5/4-\epsilon'}}
	\end{align*}
for all $\epsilon' > 0$ as desired. 
}
\fi
\end{proof}
\end{lem}
Again, Corollary \ref{corEnd} proves the conjecture stated in \cite{REU2016}.
\begin{cor} \label{corEnd}
Let $5 \leq p,q \leq M$ be randomly chosen distinct primes and let $N = pq$. Let $\Emod{N}$ be a randomly chosen elliptic curve with good reduction at $p$ and $q$ such that $(p+1-a_p), (q+1-a_q) \mid (N+1-a_N)$. Then\begin{align*}
	\lim_{M \rightarrow \infty} \Pr[a_p \text{ or } a_q \text{ is not 1}] = 0
\end{align*} 
and
\begin{align*}
	\lim_{M \rightarrow \infty} \Pr[\#\Emod{N} = N+1-a_N] = 1.
\end{align*}
\begin{proof}
In the case when $E$ is a random elliptic curve with good reduction at $p$ and $q$, not necessarily with $(p+1-a_p),(q+1-a_q) \mid (N+1-a_N)$, Lemmas \ref{LemmaLBound} and \ref{LemmaUBound} show that
	\begin{align*}
		\frac{ \Pr[a_p \text{ or } a_q \text{ is not 1 and } (p+1-a_p), (q+1-a_q) \mid N+1-a_N]}{\Pr[a_p = a_q = 1]}  
		&\ll \frac{1}{M^{1/4-2\epsilon}}.
	\end{align*}
Thus, assuming that $E$ satisfies $(p+1-a_p), (q+1-a_q) \mid (N+1-a_N)$, 
	\begin{align*}
	\lim_{M \rightarrow \infty} \Pr[a_p \text{ or } a_q \text{ is not 1}] = 0.
	\end{align*} 
Since $\#\Emod{N} \neq N+1-a_N$ implies that $a_p$ or $a_q$ is not $1$,
	\begin{align*}
		\lim_{M \rightarrow \infty} \Pr[\#\Emod{N} \neq N+1-a_N] = 0,
	\end{align*}
and so $\lim_{M \rightarrow \infty} \Pr[\#\Emod{N} = N+1-a_N] = 1$.

\end{proof}
\end{cor}

\section*{Acknowledgment}
  \noindent The authors wish to acknowledge the reviewer for the useful comments that helped to improve the paper.

\appendix
\section{Point multiplication modulo $N$} \label{SectionMult}
Let $k$ be a field and $E/k: y^2 = x^3 + Ax + B$ with $A,B \in k$ an elliptic curve. One can define the division polynomial $\psi_n = \psi_n(x,y)$ as follows:
\begin{align*}
	\psi_0 &= 0 \\
	\psi_1 &= 1 \\
	\psi_2 &= 2y \\
	\psi_3 &= 3x^4 + 6Ax^2 + 12Bx - A^2 \\
	\psi_4 &= 4y(x^6 + 5Ax^2 + 20Bx^3 - 5A^2x^2 - 4ABx - 8B^2 -A^3) \\
	\psi_{2m+1} &= \psi_{m}^3 \psi_{m+2} - \psi_{m+1}^3 \psi_{m-1},\text{ for } m \geq 2 \\
	2y\psi_{2m} &= \psi_{m}(\psi_{m+2}\psi_{m-1}^2 - \psi_{m-2} \psi_{m+1}^2),\text{ for } m \geq 3.
\end{align*}
Given a point $P = (x,y) = [x:y:1]$ on $E/k$ and a nonnegative integer $n$, the projective coordinates of $nP$ are given as
\begin{align*}
	nP = \left[\phi_n \psi_n: \omega_n: \psi_n^3 \right],
\end{align*}
where 
\begin{align*}
\phi_n = x \psi_n^2 - \psi_{n+1}\psi_{n-1} \qquad \text{and} \qquad \omega_n = \frac{\psi_{n+2}\psi_{n-1}^2 - \psi_{n-2} \psi_{n+1}^2}{4y}.
\end{align*}
With these definitions in mind, one can show Lemma \ref{LemmaDivPoly} below.
\begin{lem} \label{LemmaDivPoly}
	Let $p$ be a prime, let $n$ be an integer and let $E/\Q$ be an elliptic curve with good reduction at $p$. If $P = (x,y)$ is a point of $\Emod{p}$, then $nP = \zero$ if and only if $\psi_n(x,y) = 0$. 
	\begin{proof}
		$nP$ is $\zero$ if and only if $\psi_n^3$ is $0$, which is true if and only if $\psi_n = 0$. 
	\end{proof}
\end{lem}
One can define the division polynomials in the same way over $\Z/N\Z$ in place of $k$ and the projective coordinates of $nP$ can be computed the same way as well. The following gives an exact criterion as to when a multiple of a point of $\Emod{N}$ is $\zero$. 
\begin{prop}
	Let $N > 1$ be an integer, let $n$ be an integer and let $E/\Q$ be an elliptic curve with good reduction at every prime dividing $N$. If $P = (x,y)$ is a point of $\Emod{N}$, then $nP \equiv \zero \pmod{N}$ if and only if $\psi_n(x,y) \equiv 0 \pmod{N}$. 
	\begin{proof}
		If $\psi_n(x,y) \equiv 0 \pmod{N}$, then $\phi_n\psi_n \equiv \psi_n^3 \equiv 0 \pmod{N}$. Conversely, suppose that $nP \equiv \zero \pmod{N}$. For each prime $p$ dividing $N$, $\psi_{n+1}$ and $\psi_{n-1}$ are nonzero modulo $p$ by Lemma \ref{LemmaDivPoly} because $(n \pm 1)P \equiv \pm P \not\equiv \zero \pmod{p}$. On the other hand, $\psi_n \equiv 0 \pmod{p}$ for each prime $p$ dividing $N$ also by Lemma \ref{LemmaDivPoly}. Therefore, $\phi_n = x\psi_n^2 -\psi_{n+1}{\psi_{n-1}}$ is invertible modulo $N$. Since $nP \equiv \zero \pmod{N}$, $\phi_n\psi_n$ must be $0$ modulo $N$, so $\psi_n \equiv 0 \pmod{N}$ as desired. 
	\end{proof}
\end{prop}
When $2y$ is invertible modulo $N$, this is a convenient means to tell whether a multiple of a point $P \in \Emod{N}$ is $\zero$ and to compute the actual value of the multiple. However, $\psi_n$, for even $n \geq 6$, and $\omega_n$, for general $n$, are infeasible to compute if $2y$ is not invertible modulo $N$ as they are defined because their computations involve a division by $2y$. Fortunately, one can tweak the definitions of these polynomials to avoid inversions by $y$. Define $\hat{\psi}_n(x,y)$ as
\begin{align*}
	\hat{\psi_n}(x,y) = \begin{cases} \frac{\psi_n(x,y)}{2y} &\text{if } n \text{ is even} \\
																		\psi_n(x,y) &\text{if } n \text{ is odd}. \end{cases}
\end{align*}
Note that $\psi_0, \psi_2$ and $\psi_4$ are all multiples of $2y$ as polynomials. Moreover,
\begin{align*}
	\hat{\psi}_{2m+1} &= \begin{cases} 16y^4 \hat{\psi}_m^3 \hat{\psi}_{m+2} - \hat{\psi}_{m+1}^3 \hat{\psi}_{m-1} &\text{if } m \geq 3 \text{ is even} \\
													\hat{\psi}_m^3 \hat{\psi}_{m+2} - 16y^4 \hat{\psi}_{m+1}^3 \hat{\psi}_{m-1} &\text{if } m \geq 3 \text{ is odd}. \end{cases} \\
	\hat{\psi}_{2m} &= \hat{\psi}_m \left( \hat{\psi}_{m+2} \hat{\psi}_{m-1}^2 - \hat{\psi}_{m-2}\hat{\psi}_{m+1}^2 \right), \text{ for } m \geq 2
\end{align*}
and
\begin{align*}
\phi_n &= \begin{cases} 4xy^2 \hat{\psi}_n^2 - \hat{\psi}_{n+1} \hat{\psi}_{n-1} &\text{if } n \text{ is even} \\
											x \hat{\psi}_n^2 - 4y^2 \hat{\psi}_{n+1}\hat{\psi}_{n-1} &\text{if } n \text{ is odd} \end{cases} \\
\omega_n &= \begin{cases} \frac{1}{2} \left( \hat{\psi}_{n+2} \hat{\psi}_{n-1}^2 - \hat{\psi}_{n-2} \hat{\psi}_{n+1}^2 \right) &\text{if } n \text{ is even} \\
													y \left( \hat{\psi}_{n+2} \hat{\psi}_{n-1}^2 - \hat{\psi}_{n-2} \hat{\psi}_{n+1}^2 \right) &\text{if } n \text{ is odd}. \end{cases}
\end{align*}
From here, one can compute
\begin{align*}
	nP = [\phi_n \psi_n : \omega_n : \psi_n^3] = \begin{cases} \left[ 2y \phi_n \hat{\psi}_n, \omega_n, (2y\hat{\psi}_n)^3 \right] &\text{if } n \text{ is even} \vspace{0.2cm}\\ 
																														\left[ \phi_n \hat{\psi}_n, \omega_n, \psi_n^3 \right] &\text{if } n \text{ is odd}. \end{cases}
\end{align*}

\section{Examples} \label{SectionMuller}
Using the methods used in Appendix \ref{SectionMult}, we note some discrepancies between \cite[Table 2]{Muller} and our computational results just as in Example \ref{ExampleMuller} in \cite{Muller}. These are all claimed to be strong, but not Euler, elliptic pseudoprimes. \par
Let $N = 9090870127122419 = 61 \cdot 997 \cdot 1289 \cdot 3851 \cdot 30113$, $E: y^2 = x^3 - 5x$ and $P = (5,10)$. $N$ is not even an elliptic pseudoprime for $(E,P)$ because $(N+1)P$ does not reduce to $\zero$ modulo $997,1289,3851$ and $30113$. \par
Let $N = 32759 = 17 \cdot 41 \cdot 47$, $E: y^2 = x^3 - 3500x - 98000$ and $P = (84,448)$. \cite{Muller} states $P = (84,884)$, but this is not on $\Emod{N}$. Rather, $P$ should be understood as $(84, 448)$. Moreover, \cite{Muller} indicates that $\left( \frac{N+1}{2^3} \right) P$ is congruent to $(2345,0)$ modulo $N$, but it seems to be congruent to $(30041,29274)$ modulo $N$. Note that $29274$ is divisible by $17$ and $41$, but not $47$, so $N$ is not a strong elliptic pseudoprime for $(E,P)$.

\end{document}